\input amstex
\documentstyle{amsppt}
\magnification=1200
\hcorrection{0.2in}
\vcorrection{-0.4in}


\define\op#1{\operatorname{#1}}
\NoRunningHeads
\NoBlackBoxes
\topmatter
\title Collapsed Manifolds with local Ricci bounded Covering Geometry
\endtitle
\date
\enddate
\author Xiaochun Rong  \footnote{A part of this work was done during the author's sabbatical leave in the fall of 2018 and spring of 2022 at Capital Normal University, partially supported by NSFC 11821101, BNSF Z19003, and
a research fund from Capital Normal University. \hfill{$\,$}}
\endauthor
\address Mathematics Department, Rutgers University
New Brunswick, NJ 08903 USA
\endaddress
\email rong\@math.rutgers.edu
\endemail
\address Mathematics Department, Capital Normal University, Beijing,
P.R.C.
\endaddress
\abstract For $\rho, v>0$, we say that an $n$-manifold $M$ satisfies
local $(\rho,v)$-bound Ricci covering geometry, if Ricci curvature $\op{Ric}_M\ge -(n-1)$, and for all $x\in M$,
$\op{vol}(B_\rho(\tilde x))\ge v>0$, where $\tilde x$ is an inverse image of $x$ on the (local) Riemannian universal cover
of the $\rho$-ball at $x$.

In this paper, we extend the nilpotent fiber bundle theorem of Cheeger-Fukaya-Gromov on
a collapsed $n$-manifold $M$ of bounded sectional curvature to $M$ of a local
$(\rho,v)$-bound Ricci covering geometry, and $M$ is close to a non-collapsed Riemannian
manifold of lower dimension. The nilpotent fiber bundle theorem significantly improves fiber bundle theorem in \cite{Hu}, and it strengthens a nilpotent fiber bundle
seen from \cite{NZ} and implies the torus bundle in \cite{HW}, which are obtained under additional local or global topological conditions, respectively.

Our construction of a nilpotent fibration requires a new proof for a result in \cite{HKRX}:
if an $n$-manifold $M$ with local $(1,v)$-bound Ricci covering
geometry has diameter $<\epsilon(n,v)$, a constant depends on $n$ and $v$, then $M$ is diffeomorphic to an infra-nilmanifold. The proof in \cite{HKRX} is to show that the Ricci flows
produces an almost flat metric, thus the result follows from the Gromov's theorem on almost flat manifolds. The new proof is independent of the Gromov's theorem, thus has which as a corollary.
If the first Betti number $b_1(M)=n$, then $M$ satisfies a $(1,v)$-bound Ricci covering geometry, thus $M$ is diffeomorphic to a standard torus (\cite{Co2}).
\endabstract
\endtopmatter
\document

\head 0. Introduction
\endhead

\vskip4mm

A complete $n$-manifold $M$ is called $\epsilon$-collapsed, if every $1$-ball on $M$ has volume less than $\epsilon$, while normalizing a bound on curvature to exclude a rescaling.

In this paper, we will extend the nilpotent fiber bundle theorem of Cheeger-Fukaya-Gromov on
a collapsed $n$-manifold $M$ of bounded sectional curvature (Theorem 0.2) to $M$ of a local
$(\rho,v)$-bound Ricci covering geometry (see below), and $M$ is close, in the Gromov-Hausdorff distance,
to a non-collapsed Riemannian manifold of lower dimension (Theorem B). The nilpotent fiber bundle theorem significantly improves fiber bundle theorem in \cite{Hu} (Theorem 0.4), it
much strengthens a nilpotent fiber bundle seen from \cite{NZ} (Theorem 0.9) and implies the torus bundle in \cite{HW} (Theorem 0.9), which are obtained under additional local and global topological restriction respectively.

Our construction of a nilpotent fibration requires a new proof for a result in \cite{HKRX}:
if an $n$-manifold $M$ with local $(1,v)$-bound Ricci covering
geometry has diameter $<\epsilon(n,v)$, a constant depends on $n$ and $v$, then $M$ is diffeomorphic to an infra-nilmanifold (Theorem 0.3). The proof in \cite{HKRX} is to show that the Ricci flows
produces an almost flat metric, thus the result flows from the Gromov's theorem on almost flat manifolds. The new proof is independent of the Gromov's theorem (Theorem A), thus has which as a corollary (Corollary 0.5). Moreover, if $b_1(M)=n$ implies that $M$ satisfies $(1,v)$-bound Ricci covering geometry ((C1) in Theorem C), thus $M$ is diffeomorphic to a standard torus (Theorem 0.9), the maximal first Betti number rigidity of almost non-negative Ricci curvature (\cite{Co2}).

Let's start with some back ground. In the study of $\epsilon$-collapsed manifolds with normalizing bound on sectional curvature, $|\op{sec}_M|\le 1$, the following almost flat manifolds theorem of Gromov has been a cornerstone:

\proclaim{Theorem 0.1} {\rm (Almost flat manifolds \cite{Gr}, \cite{Ru})} Given $n>0$, there exist constants,
$\epsilon(n), w(n)>0$, such that if the scaling invariant of a compact $n$-manifold $M$,
$(\max |\op{sec}_M|) (\op{diam}(M))^2<\epsilon(n)$, or equivalently,
$$|\op{sec}_M|\le 1,\quad \op{diam}(M)<\epsilon(n),$$
then $M$ is diffeomorphic to an infra-nilmanifold, $N/\Gamma$, where $N$ is a
simply connected nilpotent Lie group, $\Gamma<N\rtimes \op{Aut}(N)$ is a discrete subgroup such that
$[\Gamma:\Gamma\cap N]\le w(n)$.
\endproclaim

In \cite{Gr}, Gromov proved that a bounded normal covering of $M$ is diffeomorphic to a nilmanifold, and that $M$
is diffeomorphic to an infra-nilmanifold was due to \cite{Ru} via constructing
a flat connection with a parallel torsion on $M$.

Starting with Theorem 0.1, Cheeger-Gromov (\cite{CG1,2}), Fukaya (\cite{Fu1-3}), and Cheeger-Fukaya-Gromov (\cite{CFG}) investigated an $\epsilon$-collapsed manifold $M$
with $|\op{sec}_M|\le 1$, and the main discovery is a compatible
local elementary nilpotent structures on $M$, every orbit has a positive dimension and points to all the collapsed directions (Theorem 1.2.2). Conversely, if a manifold admits such a compatible elementary nilpotent structures, then one constructs an one parameter family of invariant metrics $g_\epsilon$ on $M$ ($\epsilon\in (0,1]$) with $|\op{sec}_{g_\epsilon}|\le 1$ and
any orbit collapses to a point, which implies that for any $x\in M$, $\op{vol}_{g_\epsilon}(B_1(x))<\epsilon$
(Theorem 1.2.1, cf. \cite{CR}, \cite{CG1}, \cite{Fu1}).

In the case that a collapsing is taking place with a bounded diameter (see blow), the nilpotent structure is pure, which has the following alternative formulation:

\proclaim{Theorem 0.2} {\rm (\cite{Fu1,2}, \cite{CFG})} Assume a sequence of compact $n$-manifolds converging in the Gromov-Hausdorff topology, $M_i@>\op{GH}>> X$, such that
$$|\op{sec}_{M_i}|\le 1,\quad \op{diam}(M_i)\le d, \quad \dim(X)<n.$$

\noindent {\rm (0.2.1)} {\rm (Nilpotent fiber bundle)} If $X$ is a Riemannian manifold, then
$d_{\op{GH}}(M_i,X)<\epsilon(X)$ (a constant depending on $X$) implies a smooth fiber bundle map, $f_i: M_i\to X$,
with fiber an infra-nilmanifold and an affine structural group, satisfying

{\rm (0.2.1.1)} $f_i$ is an $\epsilon_i$-Gromov-Hausdorff approximation (briefly GHA), $\epsilon_i\to 0$.

{\rm (0.2.1.2)} $f_i$  is an $\epsilon_i$-Riemannian submersion i.e., for vector $x_i$ orthogonal to a $f_i$-fiber, $e^{-\epsilon_i}\le \frac{|df_i(\xi_i)|}{|\xi_i|}\le e^{\epsilon_i}$.

{\rm (0.2.1.3)} The second fundamental form of $f_i$-fibers, $|\op{II}_{f_i}|\le c(n)$.

\vskip2mm

\noindent {\rm (0.2.2)} {\rm (Singular nilpotent fibration)} If $X$ is not a Riemannian manifold, then the frame bundle of $M_i$ admits a smooth $O(n)$-invariant fiber
bundle map, $F_i\to (F(M_i),O(n))@>f_i>>(Y,O(n))$, with $Y$ a $C^{1,\alpha}$-Riemannian manifold, $F_i$ a compact nilmanifold
and an affine structure group,  such that the following diagram commutes,
$$\CD F_i \to (F(M_i),O(n))@> f_i>>(Y,O(n))\\
@VV \op{proj}_iV @VV \op{proj}V\\
\bar F_i\to M_i@>\bar f_i>>X=Y/O(n)
\endCD$$
where $f_i$ satisfies {\rm (0.2.1.1)-(0.2.1.3)}, $\bar f_i$ is the descending map of $f_i$ and every $\bar F_i=\op{proj}_i(F_i)$ is diffeomorphic to an infra-nilmanifold.
\endproclaim

For the existence of a mixed nilpotent structure on an $\epsilon$-collapsed $n$-manifold with $|\op{sec}|\le 1$ (but no bound on $\op{diam}(M)$), see Section 1.2, and for applications of the existence of nilpotent structures on collapsed manifolds with bounded sectional curvature, see \cite{Fu3}, \cite{Ro1,2}.

Let's turn to a collapsed manifold with bounded Ricci curvature. The weak curvature regularity
bring in a complexity of underlying structures. For instance, it is unlikely to expect a classification on an almost Ricci flat manifolds
(i.e., $|\op{Ric}_M|\cdot (\op{diam}(M))^2<\epsilon$); partly because which contains all compact Ricci flat manifolds (comparing to Theorem 0.1), nor to expect any fibration structure as in (0.2.1), if replacing $|\op{sec}_M|\le 1$ with $|\op{Ric}_M|\le n-1$ (\cite{An}).

In view of the above, most work in recent investigation has been focusing on situations with various additional conditions (\cite{CRX1,2}, \cite{DWY}, \cite{HKRX}, \cite{HRW}, \cite{HW1}, \cite{Ka}, \cite{NT1,2}, \cite{NZ},  \cite{PWY}, etc.), in finding a nearby metric
with bounded sectional curvature via local or global smoothing techniques (\cite{Ab}, \cite{DWY}, \cite{Ha}, \cite{HW2}, \cite{Per}, \cite{PWY}), so as to obtain a underlying nilpotent structure from the work of Cheeger-Fukaya-Gromov.

Around 2016, the author proposed to investigate the class of collapsed manifolds satisfying a local {\it $(\rho, v)$-bound Ricci covering geometry}; which includes all manifolds in the above work, as well as ones admitting no nearby metric with bounded sectional curvature (Example 1.5). We conjectured
that this class of collapsed manifolds should pose similar nilpotent structures (Conjecture 0.6). The proposal has been partially supported by recent construction of collapsing Calabi-Yau metrics (\cite{GW}, \cite{GTZ1,2}, \cite{HSVZ}, \cite{Zh}, etc), where a local $(\rho,v)$-bound Ricci covering geometry are satisfied on $M$, or a region in $M$ away from compact subset of lower dimension.

For $\rho>0$, we call the metric ball, $B_\rho(\tilde x)$, a {\it local rewinding} of $B_\rho(x)$, $x\in M$, where $\pi: (\widetilde{B_\rho(x)},\tilde x)\to (B_\rho(x),x)$ is the (incomplete) Riemannian universal cover.  We call the number, $\widetilde{\op{vol}}(B_\rho(x)):=\op{vol} (B_\rho(\tilde x))$, the local rewinding volume of $B_\rho(x)$.  Clearly, $\widetilde{\op{vol}}(B_\rho(x))=\op{vol}(B_\rho(x))$ if and only if $B_\rho(x)$ is simply connected. If $\op{vol}(B_\rho(x))$ is small and $\widetilde{\op{vol}}(B_\rho(x))$ is not small, we say that $B_\rho(x)$ is collapsed with
local bounded covering geometry. We say that $M$ satisfies local $(\rho,v)$-bound Ricci covering geometry if
$$\op{Ric}_M\ge -(n-1),\quad \widetilde{\op{vol}}(B_\rho(x))\ge v>0,\quad \forall\, x\in M,$$
or briefly, $M$ satisfies local Ricci bounded covering geometry.

First, any collapsed $n$-manifold with $|\op{sec}|\le 1$ satisfies local $(\rho,v)$-bound covering geometry with $\rho$ and $v$ depending only on $n$ (\cite{CFG}). Secondly, if
a collapsing is taking place along a nilpotent fibration i.e., each fiber collapses to
a point while its normal slice is not collapsed, then the Riemannian universal cover of a local trivialization is not collapsed. Thirdly, a collapsed manifold with local/global $(\rho,v)$-bound Ricci covering geometry may admit no nearby metric of bounded sectional curvature (Example 1.5.1), or the present smoothing techniques fail to apply (Theorem 0.4).

The main goal of our investigation is to establish a similar nilpotent structure on a collapsed manifold with
local $(\rho,v)$-bound  Ricci covering geometry; whose regular (resp. singular) part coincides (resp. is different)
with that of a nilpotent structure in \cite{CFG} ( Conjecture 0.6 below). Following
a guide in the Cheeger-Fukaya-Gromov theory on collapsed manifolds with bounded sectional curvature, the first step is to classify maximally collapsed manifold (i.e., $\op{diam}(M)$ is small) with a local Ricci bounded covering geometry
.

\proclaim{Theorem 0.3} {\rm (Maximally collapsed manifolds with Ricci bounded covering geometry, \cite{HKRX})} Given $n, v>0$, there exist constants, $\epsilon(n,v), w(n)>0$, such that if a compact $n$-manifold $M$ satisfies
$$\op{Ric}_M\ge -(n-1),\quad \op{vol}(B_1(\tilde p))\ge v,\quad \op{diam}(M)<\epsilon(n,v),$$
then $M$ is diffeomorphic to an infra-nilmanifold, $N/\Gamma$, where $\tilde p$ is a point in
the Riemannian universal cover of $M$, $N$ is a
simply connected nilpotent Lie group, $\Gamma<N\rtimes \op{Aut}(N)$ is a discrete subgroup such that
$[\Gamma:\Gamma\cap N]\le w(n)$.
\endproclaim

Note that for a maximally collapsed manifold $M$, the local Ricci bounded covering geometry is equivalent to
that the Riemannian universal cover of $M$ is not collapsed.

The proof in \cite{HKRX} showed that the Ricci flows on $M$ exists for a definite time which yields an almost flat
metric, thus Theorem 0.3 follows from Theorem 0.1. The time estimate relies on the Perel'man's pseudo-locality (\cite{Per})
on the Riemannian universal cover $\tilde M$, which requires
that every point $\tilde x\in \tilde M$,  $\op{vol}(B_1(\tilde x))$ almost equals to $\op{vol}(\b B^n_1(0))$
(\cite{CM}), and which is verified based on the Cheeger-Colding's work on a non-collapsed Ricci limit space
(\cite{CC1-3}).

The next step is to generalize Theorem 0.2 to a collapsed manifold with local $(\rho,v)$-bound  Ricci covering geometry;
which has been partially done in the case that $X$ is a Riemannian manifold (see (0.2.1)).

\proclaim{Theorem 0.4} {\rm (Fiber bundles, \cite{Hu})} Given $n, v>0$, there exists a constant $\epsilon(n,v)>0$ such that if
compact $n$-manifold $M$ and $m$-manifold $N$ satisfy, $\forall\, x\in M$, $0<\epsilon\le \epsilon(n,v)$,
$$\op{Ric}_M\ge -(n-1),\, \widetilde{\op{vol}}(B_1(x))\ge v, \, |\op{sec}_N|\le 1,\, \op{injrad}(N)\ge 1, \,
d_{\op{GH}}(M,N)<\epsilon,$$
then there is a smooth fiber bundle map, $f: M\to N$, which is also an $\Psi(\epsilon|n,v)$-GHA, where $\Psi(\epsilon|n,v)\to 0$ as $\epsilon\to 0$ while $n$ and $v$ are fixed.
\endproclaim

Note that smoothing methods mentioned in the above work unlikely apply under the conditions of Theorem 0.4. In the proof (\cite{Hu}), $f$ is constructed by gluing locally defined $(m,\delta)$-splitting maps (\cite{CC1}) via the standard center of mass technique associate to a partition of unity on $N$, and the non-degeneracy of $f$ relies on the non-degeneracy
of any $(n,\delta)$-splitting map (\cite{CJN}).

Comparing with (0.2.1), the significant structural information missing in Theorem 0.4 is that a $f$-fiber should be diffeomorphic to an infra-nilmanifold, and
the structural group should be reduced to an affine group; the two properties are expected and have been indispensable
in applications (\cite{Ro1}). Indeed, in the case that the Riemannian universal cover $\tilde M_i$ is not collapsed, a singular nilpotent fibration with an affine structural group was obtained in \cite{HKRX} by a smoothing method via Ricci flows, similar to the proof of Theorem 0.3.

In contrasting to that (0.2.1) follows from Theorem 0.1, Theorem 0.3 fails to apply to $f$-fiber
in Theorem 0.4 because the intrinsic metric on a $f$-fiber unlikely satisfies Theorem 0.3.
Hence, to conclude that a fiber in Theorem 0.4 is diffeomorphic to an infra-nilmanifold and
there is an affine structural group, one needs to find a new proof of Theorem 0.3 without
Theorem 0.1.

We now begin to state main results in this paper.


\proclaim{Theorem A} Theorem 0.3 can be proved independent of Theorem 0.1.
\endproclaim

A significance of Theorem A is that the proof naturally extends to a proof of a nilpotent fibration theorem (Theorem B). A bi-product from the proof of Theorem A is that one explicitly constructs a nearby left invariant almost flat metric (see completion of proof
of Theorem 0.3').

A consequence of Theorem 0.1 is that the Riemannian universal covering of $M$ of $|\op{sec}_M|\le 1$ and $\op{diam}(M)<\epsilon(n)$ is not collapsed. An elementary proof independent of Theorem 0.1 was recently given (see (8.4) in \cite{Ro3}). Hence, a consequence of Theorem A is

\proclaim{Corollary 0.5} Theorem 0.3 implies Theorem 0.1.
\endproclaim

Our approach to Theorem A is partially motivated by a new proof for Theorem 0.1 in \cite{Ro4}, which is by induction on dimension based on a recent criterion of nilmanifolds: a compact manifold $M$ is a nilmanifold if and only if $M$ admits an iterated principal circle bundles (\cite{Na}, \cite{Be}). Because the inductive proof relies on smoothing techniques we are not allowed to use, we formulate an alternative criterion for our purpose (Theorem 2.1), which
can be verified using tools available on manifolds of Ricci curvature bounded below (a brief outline is given at the end of introduction).

\proclaim{Theorem B} {\rm (Nilpotent fiber bundles)} Let the assumptions be as in Theorem 0.4. Then there is a
smooth bundle map, $f: M\to N$, such that $f$ is an $\Psi(\epsilon|n,v)$-GHA, a $f$-fiber is diffeomorphic
to an infra-nilmanifold, and the structural group is affine.
\endproclaim

Theorem B substantially improves Theorem 0.4 by asserting a fiber an infra-nilmanifold and an affine structural group.

We will present two proofs of Theorem B: one is to run the proof of Theorem A to a fiber
(in a tubular neighborhood) in Theorem 0.4, and the other is to, by locally running
the proof of Theorem A, construct a locally finite open cover for $M$, each supports an elementary N-structure (see Section 1.3), for which to apply gluing techniques
in \cite{CFG} to patch these elementary N-structures to a desired nilpotent fibration.

Theorem B extends (0.2.1) to collapsed manifolds with local $(\rho,v)$-bound Ricci covering geometry. Concerning a possible extension of (0.2.2), we recall the following conjecture (\cite{Ro4}, \cite{HRW}):

\example{Conjecture 0.6} Let a sequence of compact $n$-manifolds, $M_i@>\op{GH}>> X$ (compact), with $\op{Ric}_{M_i}\ge -(n-1)$ and $\op{vol}(M_i)\to 0$.

\noindent (0.6.1) {\rm(Singular nilpotent fibration)} Assume all $x_i\in M_i$ are local rewinding $(\rho,\epsilon)$-Reifenberg points i.e.,
$d_{\op{GH}}(B_{i,s}(\tilde x_i),\b B^n_s(0))<s\epsilon$, $0<s\le \rho$, where $(\widetilde{B_\rho(x_i)},\tilde x_i)\to (B_\rho(x_i),x_i)$ is
a (local) Riemannian universal covering map. Then for $i$ large, there
is a singular nilpotent fibration, $f_i: M_i\to X$, which satisfies the commutative diagram in (0.2.2).

\noindent (0.6.2)  {\rm(Generalized singular nilpotent fibration)} Assume $M_i$ satisfies local $(\rho,v)$-bound  Ricci covering geometry.
Then for $i$ large, there is a continuous map, $f_i: M_i\to X$, $f_i$ is an $\epsilon_i$-GHA, and for any $x\in X$, $f^{-1}_i(x)$ is
finite quotient of an infra-nilmanifold, and at $x_i\in f^{-1}_i(x)$, any short geodesic loop is shortly homotopic to a nontrivial element
in $\pi_1(f^{-1}_i(x))$.
\endexample

Note that in the case that $X$ is a Riemannian manifold, the conditions in (0.6.1) and (0.6.2) are equivalent,
and Conjecture 0.6 holds in this case (Theorem B).  

Roughly speaking, (0.6.2) says that the regular part of the singular nilpotent fibration coincides with one in
Theorem 0.2, while the singular part is slightly `singular' than one in Theorem 0.2 because a singular fiber
is an infra-nilorbifold.

Next, we will give applications of Theorem B relating to recent work in \cite{NZ} and \cite{HW1,2} on a collapsed manifold $M$ with Ricci curvature bounded in two sides or below, satisfying a local or global topological condition, respectively. In each work,  mainly using analytic tools, it was proved that $M$ satisfies a local or global bounded covering geometry, respectively. We will extend such bound covering property in a more general setting
(Theorem C), and we will present a proof based on elementary tools in metric geometry.

The generalized Margulis lemma in \cite{KW} asserts that there exist constants, $\epsilon(n), c(n)>0$, such that for any point $x$
in a complete $n$-manifold $M$ of $\op{Ric}_M\ge -(n-1)$, the subgroup of $\pi_1(B_1(x),x)$, $\Gamma_x(\epsilon)$, generated
by loops at $x$ of length $<\epsilon\le \epsilon(n)$ contains a nilpotent subgroup of index $\le c(n)$ (Theorem 1.4.2).  In particular,
$\Gamma_x(\epsilon)$ contains a torsion-free nilpotent subgroup of finite index (\cite{Ra}).  The rank of
a torsion-free nilpotent group is equal to the minimal number of generators, and the rank
of $\Gamma_x(\epsilon)$ is defined as the rank of the torsion-free nilpotent subgroup of finite index.
For the notions on pseudo-groups and their action, see Subsection 1.1 and Remark 1.1.5.

Let $X$ be a Ricci limit space. According to \cite{CoN}, there is a unique integer $k$ such that the subset of
$X$ consisting of points whose tangent cone is unique and is isomorphic to $\Bbb R^k$ has a full renormalized measure. We will use the notation that $\dim(X)=k$.

\proclaim{Theorem C} {\rm (Local/global Ricci bounded covering geometry)} Let a sequence of compact
$n$-manifolds, $(M_i,p_i)@>\op{GH}>>(X,p)$, with $\op{Ric}_{M_i}\ge -(n-1)$,  $\op{diam}(M_i)\le d$.

\vskip1mm

\noindent {\rm (C1) (Local rank estimate)} For $x_i\in M_i$, assume the following commutative diagram:
$$\CD (\bar B_1(\tilde x_i),\Gamma_{x_i}(\epsilon))@>\op{eqGH}>>(\bar B_1(\tilde x),G)\\
@V \pi_i|_{B_1(\tilde x_i)}VV @V \op{proj}VV\\
\bar B_1(x_i)@>\op{GH}>>B_1(x)=\bar B_1(\tilde x)/G,\endCD$$
where $\pi_i: (\widetilde {B_2(x_i)},\tilde x_i)\to (B_2(x_i),x_i)$ is a Riemannian universal cover map,
$\Gamma_{x_i}(\epsilon)$ acts pseudolly on $B_1(\tilde x_i)\subset \widetilde{B_2(x_i)}$, and $G$, the limit of $\Gamma_{x_i}(\epsilon)$,
is a pseudo-group of isometries.
Then $\op{rank}(\Gamma_{x_i}(\epsilon))\le \dim(G)$, and $\op{rank}(\Gamma_{x_i}
(\epsilon))=n-\dim(X)$ implies that $\op{vol}(B_1(\tilde x_i))\ge v(X)>0$.

\vskip1mm

\noindent {\rm (C2) (Global abelian rank estimate)} Given the following commutative diagram,
$$\CD (\tilde M_i,\tilde p_i,\Gamma_i)@>\op{eqGH}>>(\tilde X,\tilde p,G)\\
@V \pi_iVV @V \op{proj}VV\\
(M_i,p_i)@>\op{GH}>>(X=\tilde X/G,p),\endCD$$
the first Betti numbers, $b_1(M_i)-b_1(X)\le \dim(G)$, and
$b_1(M_i)-b_1(X)=n-\dim(X)$ implies that $\op{vol}(B_1(\tilde p_i))\ge v(X)>0$.
\endproclaim

In (C1), the pseudo-group $G$ contains a neighborhood of the identity that isomorphically embeds into a Lie group $\hat G$, and $\dim(G):=\dim(\hat G)$.

\remark{Remark \rm 0.7}  \noindent (0.7.1) If one assumes that $X$ is Riemannian manifold i.e.,
a collapsed $M$ is $d_{GH}$-close to a non-collapsed Riemannian manifold, then one may state
the non-collapsing on a local or global Riemannian universal cover in (C1) or (C2),
without involving a GH-convergent sequence (\cite{NZ}, \cite{HW1}).

(0.7.2) If $X$ is a compact Riemannian $m$-manifold, the case of (C1) that
$|\op{Ric}_M|\le (n-1)$ was obtained in \cite{NZ}, and if $X$ is an orbifold, (C2) was obtained
in \cite{HW2}. 
\endremark

A consequence of Theorem A and (C1) is

\proclaim{Therem 0.8} {\rm(\cite{Co2})} Given $n$, there exists $\delta(n)>0$ such that if a compact
$n$-manifold satisfies
$$\op{Ric}_M\ge -(n-1),\quad \op{diam}(M)\le \delta(n),\quad b_1(M)=n,$$
then $M$ is diffeomorphic to a flat torus.
\endproclaim

A consequence of Theorem B and (C1) is

\proclaim{Theorem 0.9} {\rm (Local maximal rank and nilpotent fiber bundles)} Assume that compact $n$-manifold $M$ and
$m$-manifold $N$ $(m<n$) satisfy
$$\op{Ric}_M\ge -(n-1),\quad  |\op{sec}_N|\le 1,\quad \op{injrad}(N)\ge 1,\quad d_{\op{GH}}(M,N)<\epsilon\le \epsilon(n),$$
where $\epsilon(n)$ is the constant in Margulis lemma in \cite{KW}. Then
$\op{rank}(\Gamma_{x}(\epsilon))\le n-m$, and ``$=$''  for all $x\in M$ implies a smooth fiber bundle, $f: M\to N$,
with fiber diffeomorphic to an infra-nilmanifold,  and $f$ is an $\Psi(\epsilon|n)$-GHA.
\endproclaim

The case of Theorem 0.9 that $|\op{Ric}_M|\le n-1$ can be easily seen from \cite{NZ}.

A consequence of Theorem B and (C2) is:

\proclaim{Theorem 0.10} {\rm (Generalized maximal first Betti number rigidity, \cite{HW1})} Assume that compact $n$-manifold $M$ and
$m$-manifold $N$ $(m<n$) satisfy
$$\op{Ric}_M\ge -(n-1),\quad  |\op{sec}_N|\le 1,\quad \op{injrad}(N)\ge 1,\quad d_{\op{GH}}(M,N)<\epsilon\le \epsilon(n).$$
Then the first Betti numbers, $b_1(M)-b_1(N)\le n-m$, and ``$=$'' implies a smooth torus fiber
bundle, $f: M\to N$, with an affine structural group, and $f$ is an $\Psi(\epsilon|n)$-GHA.
\endproclaim

Theorem 0.10 may be viewed as a bundle version of Theorem 0.8.

Observe that the case of Theorem 0.10 with bounded sectional curvature follows easily from Theorem 0.2. A key in
the proof is to show that $\tilde M$ is not collapsed, based on which a smoothing via Ricci flows applies (\cite{HW2}).

\vskip2mm

We will conclude the introduction with an indication of our approach to Theorems A and B (see Section 2 for a detailed outline).

By the Gromov's compactness, it is equivalent to consider a sequence of compact $n$-manifolds, $M_i@>\op{GH}>>\op{pt}$, with local $(\rho,v)$-bound  Ricci covering geometry, and show that for $i$ large, $M_i$ is diffeomorphic to an infra-nilmanifold (Theorem 0.3').

The starting point in our proof is the generalized Margulis lemma in \cite{KW} (Theorem 1.4.2), by which $M_i$ has a bounded normal covering space, $M_{i,0}$, with $\pi_1(M_{i,0})$ a nilpotent group. The proof consists of two steps: (1) (the main work) we prove that $M_{i,0}$ is diffeomorphic to a nilmanifold, by exploring a sequence of sub-bundle structures on $M_{i,0}$
such that the adjacent orbit space is a principal tori bundle over the other (Theorems 2.1, 2.2). (2) we prove that $M_i$ is diffeomorphic to an infra-nilmanifold: viewing $\tilde M_i$ as a simply connected nilpotent Lie group, by the Malc\'ev rigidity there is a $\Gamma_i=\pi_1(M_i)$-action on
$\tilde M_i$ which preserves the sub-bundle structures. We will construct a $\Gamma_i$-conjugation
on $\tilde M_i$ via defining a left invariant metric on $\tilde M_i$, based on the geometry of a short basis in $\pi_1(M_{i,0})$, and applying a generalized stability of isometric Lie group actions (Lemma 3.5.1).

The sequence of sub-bundle structures on $M_{i,0}$ is constructed through a successively blowing up at every point in $M_{i,0}$, with rates obtained from a short basis at the point, and the nilpotency of $\pi_1(M_{i,0}$ guarantees that at each blow up $M_{i,0}$ is locally close (in the Gromov-Hausdorff distance) to a metric product
of a flat torus and an Euclidean space (Theorem 2.4), which is used to construct
a locally finite open cover on $M_{i,0}$, each carries an elementary T-structure. We will glue
these elementary T-structures to a fiber bundle global (\cite{CG2}, \cite{GK}, \cite{Pa1}),
and the metric regularity required for the patching are that all points on the Riemannian universal cover are uniform (arbitrarily) Reifenberg (Lemma 2.3), which is achieved based on the Cheeger-Colding-Naber theory on the degeneration of metrics with Ricci curvature bounded below  (\cite{CC1-3}, \cite{Co1,2}, \cite{CN}, \cite{CJN}).


\vskip2mm

In the proof of Theorem B, by performing a similar successive blow ups as in the proof of
Theorem A, we obtain a locally finite open cover on $M$ with elementary N-structure, $f_\alpha: U_\alpha\to B_\alpha$, where $f_\alpha$ is a smooth trivial bundle map and $\Psi(\epsilon | n,v)$-GHA, with fiber an infra-nilmanifold, such that on overlaps, elementary N-structures satisfy the patching criterion in \cite{CFG}. Hence we are able glue elementary N-structures into a nilpotent fiber bundle with an affine structural group.

Finally, we pointed out the our proofs of Theorems A and B are mainly local arguments,
hence they can be applied to a region in a collapsed manifold where a local Ricci bounded covering geometry is satisfied (e.g., \cite{GW}, \cite{GTZ1,2}, \cite{HSVZ}, \cite{Zh}).

\vskip6mm

The rest of the paper will be organized as follows:

In Section 1, we will review basic notions and results that will be used in the rest of the paper.

In Section 2, we will outline a proof of Theorem 0.3 that does not rely on Theorem 0.1.

In Section 3, we will prove Theorem A.

In Section 4, we will prove Theorem B.

In Section 5,  we will prove Theorem C, Theorem 0.8 and Theorem 0.10.

\vskip6mm

\head 1. Preliminaries
\endhead

\vskip4mm

We will review basic notions and properties that will be used in the rest of the paper. We will also present examples
mentioned in the introduction.

\vskip4mm

\subhead 1.1. Equivariant Gromov-Hausdorff convergence
\endsubhead

\vskip4mm

The reference of this subsection is [FY] (cf. [Ro2]).

Let $X_i$ and $X$ be compact length spaces. We say that $X_i$ converge to $X$ in the Gromov-Hausdorff topology,
denoted by $X_i@>\op{GH}>>X$, if there are a sequence $\epsilon_i\to 0$ and a sequence of maps,
$h_i: X_i\to X$, such that $|d_X(h_i(x_i),h_i(x_i'))-d_{X_i}(x_i,x_i)|
<\epsilon_i$ ($\epsilon_i$-isometry), and for any $x\in X$, there is $x_i
\in X_i$ such that $d_X(h_i(x_i),x)<\epsilon_i$ ($\epsilon_i$-onto), and
$h_i$ is called an $\epsilon_i$-Gromov-Hausdorff approximation, briefly,
$\epsilon_i$-GHA.

Assume that $X_i$ admits a closed group $\Gamma_i$-action by isometries. Then
$(X_i,\Gamma_i)@>\op{GH}>>(X,\Gamma)$ means that there are a closed group of
isometries $\Gamma$ on $X$, a sequence $\epsilon_i\to 0$ and
a sequence of triple of $\epsilon_i$-GHA, $(h_i,\phi_i,\psi_i)$, $h_i: X_i\to X$, $\phi_i:
\Gamma_i\to \Gamma$ and $\psi_i: \Gamma\to \Gamma_i$, such that for all $x_i\in X_i, \gamma_i\in \Gamma_i$ and $\gamma\in \Gamma$,
$$d_X(h_i(x_i), \phi_i(\gamma_i)h_i(\gamma_i^{-1}(x_i)))<\epsilon_i,\quad
d_X(h_i(x_i),\gamma^{-1}(h_i(\psi_i(\gamma)(x_i)))<\epsilon_i, \tag 1.1.1$$
where $\Gamma_i$ and $\Gamma$ are equipped with the induced metrics
from $X_i$ and $X$. We call $(h_i,\phi_i,\psi_i)$ an $\epsilon_i$-equivariant GHA.

When $X$ is not compact, then the above notion of equivariant convergence
naturally extends to a pointed version $(h_i,\phi_i,\psi_i)$: $h_i:
B_{\epsilon_i^{-1}}(p_i)\to B_{\epsilon_i^{-1}+\epsilon_i}(p)$, $h_i(p_i)=p$,
$\phi_i: \Gamma_i(\epsilon_i^{-1})\to \Gamma(\epsilon_i^{-1}+\epsilon_i)$,
$\phi_i(e_i)=e$, $\psi_i: \Gamma(\epsilon_i^{-1})\to \Gamma_i(\epsilon_i^{-1}
+\epsilon_i)$, $\psi_i(e)=e_i$, and (1.1.1) holds whenever the multiplications
stay in the domain of $h_i$, where $\Gamma_i(R)=\{\gamma_i\in\Gamma_i,\,
\,d_{X_i}(p_i,\gamma_i(p_i))\le R\}$.

\proclaim{Lemma 1.1.2}  Let $(X_i,p_i)@>\op{GH}>>(X,p)$, where $X_i$ is a  complete
locally compact length space. Assume that $\Gamma_i$ is a closed group of isometries
on $X_i$. Then there is a closed group $G$ of isometries on $X$
such that passing to a subsequence, $(X_i,p_i,\Gamma_i)@>\op{GH}>>(X,p,G)$.
\endproclaim

\proclaim{Lemma 1.1.3}  Let $(X_i, p_i, \Gamma_i)@>\op{GH}>>(X,p,\Gamma)$, where $X_i$ is
a complete locally compact length space and $\Gamma_i$ is a closed subgroup
of isometries. Then $(X_i/\Gamma_i,\bar p_i)@>\op{GH}>>(X/\Gamma,\bar p)$.
\endproclaim

Assume that the universal covering space, $\pi_i: (\tilde X_i,
\tilde p_i)\to (X_i,p_i)$, exists, and let $\Gamma_i=\pi_1(X_i,p_i)$
denote the fundamental group.

\proclaim{Lemma 1.1.4}  Let $X_i@>\op{GH}>>X$ be a sequence of compact
length metric space. Then passing to a subsequence such that the following
diagram commutes,
$$\CD (\tilde X_i,\tilde p_i,\Gamma_i)@>\op{GH}>>(\tilde X,\tilde p,\Gamma)
\\ @VV \pi_i V   @VV \pi  V \\
(X_i,p_i)@>\op{GH}>> (X,p).\endCD$$
If $X$ is compact and $\Gamma/\Gamma_0$ is discrete, then there is $\epsilon>0$
such that the subgroup, $\Gamma_i(\epsilon)$, generated by elements
with displacement bounded above by $\epsilon$ on $B_{2d}(\tilde p_i)$,
is normal and for $i$ large, $\Gamma_i/\Gamma_i(\epsilon)\overset
{\op{isom}}\to \cong \Gamma/\Gamma_0$.
\endproclaim

\remark{Remark \rm 1.1.5} Inspecting proofs of Lemma 1.1.2-1.1.4, the whole group structure of
$\Gamma_i$ is not needed; the multiplication of elements in $\Gamma_i(\epsilon_i^{-1})$
is taken into account only if the product moves a point in $B_{\epsilon_i^{-1}}(p_i)$ and remained
in the ball.   We will call $\Gamma_i(R)=\{\gamma_i\in \Gamma_i,\, d(\gamma_i(p_i),p_i)\le R\}$ is a pseudo-group, and the
$\Gamma_i(R)$-action (resp. the $\Gamma_i$-action) on $B_{\epsilon_i^{-1}}(p_i)$
a pseudo-action. Clearly, Lemma 1.1.2-1.1.4 hold when replacing $\Gamma_i$ by $\Gamma_i
(\epsilon_i^{-1})$ and the $\Gamma_i$-actions by pseudo-actions.
\endremark

\vskip4mm

\subhead  1.2. Patching atlas of elementary F-structures and N-structures
\endsubhead

\vskip4mm

In our proof of Theorem A, a basic tool is the patching technique in \cite{CG2} on gluing
an atlas of elementary F-structures on a manifold $M$ to form a (mixed) F-structure (see below),
and in the proof of Theorem B, we employ the generalized patching techniques in \cite{CFG} to glue an atlas of elementary N-structures to obtain the desired nilpotent fiber bundle.

Indeed, we need only the simplest form of the above patching technique: patching local elementary $T^k$-structure by free $T^k$-actions ($k$ is fixed) (resp. local nilpotent fiber
bundles), which are $C^1$-close on overlaps, to form an affine $T^k$-bundle i.e.,
$T^k$-fibers are orbits of local free $T^k$-actions which are commutative on the overlap (resp.
to form a nilpotent fiber bundle on $M$). Hence the structure group of the bundle is affine.

Let $M$ be a Riemannian $n$-manifold, and let $U\subseteq M$ be an open subset.  An elementary N-structure on $U$, denoted by
$(U, \Cal N)$, is defined by an $O(n)$-invariant fiber bundle on the frame bundle over $U$, $F\to (F(U),O(n))@>f>>(Y,O(n))$,
where $F=N/\Lambda$ is a compact nilmanifold, and the $O(n)$-action preserves $F$-fibers as well as
the nilpotent structure on $F$. The $O(n)$-invariance of $f$ descends the $F$-bundle to a map, $\bar f: U\to X=Y/O(n)$,
such that the following diagram commutes,
$$\CD F\to (F(U),O(n))@>f>>(Y,O(n))\\
@V\pi VV@V\op{proj}VV\\
\bar F\to U@>\bar f>>X=Y/O(n).
\endCD$$
We will call the $\pi$-image of $F$ an orbit on $U$, and $U$ decomposes into orbits.  Note that the property
that $U$ admits a pure N-structure does not dependent on the metric on $M$, because given $M$ another
metric, its frame bundle on $U$ is equivalent to $(F(U),O(n))$. If $F=T^k$, then we will call this elementary
N-structure an elementary F-structure, denoted by $(U,\Cal F)$. If the $T^k$-bundle is also principal, we will
call it an elementary T-structure; note that the free $T^k$-action on  $F(U)$ descends to a $T^k$-action on
$U$ (which may not be free).

An N-structure on a manifold $M$ is given by a collection of elementary N-structures, $(M,\Cal N)=\{(U_i,\Cal N_i)\}$,
where $\{U_i\}$ form a locally finite open cover for $M$, which satisfy the following compatibility condition:
if $U_i\cap U_j\ne \emptyset$, then $U_i\cap U_j$ is both $\Cal N_i$- and $\Cal N_j$-invariant,
and on $F(U_i\cap U_j)$, the $F_i$-bundle is a sub-bundle of the $F_j$-bundle, or vice versa. When $U_i\cap
U_j\ne \emptyset$, the $F_i$-bundle and $F_j$-bundle coincide, we will call $\Cal N$ a pure N-structure,
and otherwise $\Cal N$ is called a mixed N-structure.

For any $x\in M$, we will call the orbit at $x$ of any elementary N-structure $\Cal N_i$ containing $x$ with the maximal
dimension the $N$-orbit at $x$, denoted by $\Cal O_x$. Clearly, $M$ decomposes into disjoint N-orbits.
The minimal dimension of N-orbits is called the rank of $\Cal N$. A metric on
$M$ is called invariant, if for a canonical metric on $F(U_i)$, the induced metric on $F_i$ is left invariant.
If $M$ admits an N-structure, then $M$ always admits an invariant metric (\cite{CR}, \cite{Ro2}).

Note that the case that a pure $T^k$-structure (T-structure) with every elementary $T^k$-structure defined
by a free $T^k$-action form an affine $T^k$-bundle on $M$.

\proclaim{Theorem 1.2.1} Let $M$ be a manifold.

\noindent {\rm(1.2.1.1) (\cite{CG1})} If $M$ admits an F-structure of positive rank, then $M$ admits one parameter
family of invariant metrics, $g_\epsilon$, such that $|\op{sec}_{g_\epsilon}|\le 1$, $\op{diam}(\Cal O_x)\to 0$ as
$\epsilon\to 0$, for all $x\in M$ (we refer such a collapse as long the F-structure).

\noindent {\rm (1.2.1.2) (\cite{CR})} If $M$ admits an N-structure of positive rank, then $M$ admits one parameter
family of invariant metrics, $g_\epsilon$, such that $|\op{sec}_{g_\epsilon}|\le 1$ and $\op{diam}(\Cal O_x)\to 0$ as
$\epsilon\to 0$,  for all $x\in M$.
\endproclaim

\proclaim{Theorem 1.2.2} There  exists a constant $\epsilon(n)>0$ such that if $M$ is a complete $n$-manifold
satisfying $|\op{sec}_M|\le 1$ and $\op{vol}(B_1(x))<\epsilon\le \epsilon(n)$ for all $x\in M$, then

\noindent {\rm (1.2.2.1) (\cite{CFG})}  $M$ admits a positive rank $N$-structure, $\Cal N$, such that each
$\Cal N$-orbit has diameter $<\Psi(\epsilon)$, and there is a $C^1$  $\epsilon$-close invariant metric, $g_\epsilon$,
with higher regularity.

\noindent {\rm (1.2.2.2) (\cite{CG2})}  The center part of the nilpotent structure of $\Cal N$ form an F-structure of
positive rank.
\endproclaim

A proof of Theorem 1.2.2 is based on Theorem 0.2 and the patching techniques in \cite{CG2};  which is an extension of the stability of compact Lie group actions (\cite{Pa}, \cite{GK}).

To construct a nilpotent structure on a collapsed $M$ in Theorem 1.2.2, one  first covers $M$ with locally bounded
number of unit balls, and for each ball,
apply a local version of Theorem 0.2 on an open subset, $B_1(x_i)\subset U_i\subset B_2(x_i)$, to obtain
an elementary $N$-structure on $U_i$. Then one gets an open cover for $M$ with elementary N-structure, $\{(U_i,\Cal N_i)\}$,
and on $U_i\cap U_j\neq\emptyset$, these elementary N-structures are $C^1$-close in a suitable sense similar to that of a compact Lie group actions considered by \cite{GK} and \cite{Pa}.  Extending the patching technique in \cite{CG2},
one can glue these elementary structures to form an N-structure (\cite{CFG}).

In our proof the proof of Theorem 0.3, we perform blow-up operations on $M$ a finite number times, and each time we blow up blow up at any $x\in M$ with a scale from lengths a short basis at $x$. By assuming that $\pi_1(M)$ is nilpotent (\cite{KW}), each
blow-up operation results in a manifold with a local finite open cover by elementary $T^k$-structures defined by free almost isometric $T^k$-actions ($k$ is independent of points).

For our purpose, we will only state the gluing result (Lemmas 1.4 and 1.5 in \cite{CG2}) for the case of constructing a $T$-structure
on $M$ by patching a collection of elementary $T$-structures which are $C^1$-close on overlaps.

\proclaim{Theorem 1.2.3} {\rm(Gluing, \cite{CG2})}  Let $M$ be an $n$-manifold.  Given
$0<\epsilon<2\rho$, there exists a constant $\delta=\delta(n,\rho,\epsilon)>0$ such that $M$ admits a $T^k$-bundle
with affine structural group, if $M$ admits a collection of elementary $T^k$-structures, $\{(U_i,g_i,T^k_i)\}$ satisfying

\noindent {\rm (1.2.3.1)}  $(U_i,T^k_i)$ is an elementary T-structure with a free $T^k_i$-action and
invariant metric $g_i$ such that $\{U^\rho_i\}$ is locally finite open cover for $M$ and
$$|\op{sec}_{g_i}|\le 1, \quad \op{injrad}(g_i)\ge \frac 1{10},$$
where $U_i^\rho\subset U_i$ consisting of points distance $\ge \rho$ from $\partial \bar U_i$.

\noindent {\rm (1.2.3.2)} If $U_i^\rho\cap U_j^\rho\ne \emptyset$, then the $T^k_i$- and $T^k_j$-actions on
$U_i^\rho\cap U_j^\rho$ is $\delta$ $C^1$-close i.e.,
$$\quad e^{-\epsilon}g_i\le g_j\le e^\epsilon g_i,\quad \eta(t)=\max_{x\in U_i^\rho\cap U_j^\rho}\{d(t_i^{-1}t_j(x),x),\, t\in B_\rho(e)\}<\delta,$$
where $T^k$ the product of unit circle and $B_\rho(r)$ the $\rho$-ball of centered at the identity element $e$.
\endproclaim

\remark{Remark \rm 1.2.4} The $C^1$ closeness condition in Theorem 1.2.3 is a modification of one in \cite{GK} and \cite{Pa} where the $C^1$-closeness is for two actions of a compact Lie group on a manifold, here $U_i^\rho\cap U_j^\rho\ne \emptyset$
may neither $T^k_i$- nor $T^k_j$-invariant.  The $C^1$-closeness condition  can be achieved on a collapsed manifold $M$
in Theorem 1.2.2 via locally averaging the metric under the local almost isometric tori actions. Note that Theorem 1.2.3 can
be applied on a collapsed $M$ with Ricci curvature bounded below as long as $U_i$ and $U_j$ admit $T^k$-invariant metrics
satisfying conditions in Theorem 1.2.3, even the invariant metrics may not be $C^1$-close the induced metrics from $M$.
\endremark

A prototype of our proof of Theorem 0.3' is seen by applying Theorem 1.2.3 to a sequence, $M_i@>\op{GH}>>N$,
where $|\op{sec}_{M_i}|\le 1$ and $N$ is a compact Riemannian manifold. By Theorem 0.2, there is a fiber bundle map,
$f_i: M_i\to N$, with fiber an infra-nilmanifold and affine structural group.

Let's apply the construction of an F-structure on $M_i$ in \cite{CG2}: $M_i$ caries an F-structure $\Cal F$
of positive rank whose an orbit points the most collapsed directions. By the diameter bound (i.e., collapsing
is taking place with a bounded diameter), $\Cal F$ is pure, and because $N$ is a manifold, $\Cal F$ is
polarized and is a fiber bundle with a flat fiber. Indeed, as pointed in \cite{CFG} (p.330) that $\Cal F$ is
a part in the centers of the infra-nilmanifold fibers, thus form a sub-bundle with
fiber a flat manifold.

Inspecting the construction of F-structure in \cite{CG2}, one first gets local F-structure around each $x_i\in M_i$,
by blowing up $M_i$ by the inverse of $\ell(x_i)=\op{injrad}(M_i,x_i)$, and a limit,
$$(\ell(x_i)^{-1}M_i,x_i)@>\op{GH}>> (Y,x),$$
is a non-compact flat manifold with a compact soul, a flat manifold $F$; and an open subset at $x_i$ that is
diffeomorphic to a neighborhood of $F$ in $Y$ defines an elementary F-structure at $x_i$.  By (0.2.1.3), the induced
metrics on $f_i$-fibers are comparable, and in particular, one obtains a local finite open cover of $M_i$
by $\{U_{i,j}\simeq T^{k_s}\times \b B_{100}^{n-k_s}(0)\}$,
and on overlaps the conditions Theorem 1.2.3 is satisfied, where $\b B_r^m(0)$ denotes an $m$-dimensional
Euclidean $r$-ball.

\proclaim{Corollary 1.2.5}  Let $f_i: M_i\to N$ be as in the above. Then there is a sub-bundle,
$f_{i,s}: M_i\to B_{i,s}$, such that a $f_{i,s}$-fiber is a flat manifold and affine structural group.
\endproclaim

In \cite{CFG} (Part III, sections 5-8), Theorem 1.2.3 is generalized to glue elementary N-structures, which is more technique. Note that a gluing of elementary N-structures is carried out on local frame bundles (so one has to make gluing
$O(n)$-invariant), and the stability of a compact Lie group actions (\cite{GK}, \cite{Pa}) cannot be applied because a (non-abelian) nilpotent Lie group is not compact. Instead, the stability
of elementary N-structures also uses the Malcev's rigidity theorem for discrete co-compact
subgroups of nilpotent Lie groups.

Given two elementary N-structures, $(U_1,\Cal N_1)$ and $(U_2,\Cal N_2)$, with $U_1\cap U_2\ne \emptyset$.
Assume that there is an open ball $B_r\subset U_1\cap U_2\ne\emptyset$ such that the $\Cal N_i$-orbits on
$B_r$ is contained in $U_1\cap U_2$.  Then the two elementary N-structures are $C^1$-close, if $U_i$ admits
an invariant Riemannian metric $g_i$ with $|\op{sec}_{g_i}|\le 1$ such that
the pullback metric $\tilde g_1$ and $\tilde g_2$ on the universal cover, $\widetilde{U_1\cap U_2}$, are
$C^1$-close restricting to any ball always from the boundary of $\widetilde{U_1\cap U_2}$.

\remark{Remark \rm 1.2.6}  The above $C^1$-closeness is a simplest case for the gluing criterion in \cite{CFG} (Part III, sections 5-8). Indeed, in our proof Theorem B, the orbits on $U_i$ form a trivial fiber bundle over a unit Euclidean ball i.e., $f_i: U_i\to \b B^m_1(0)$, thus $\Gamma=\pi_1(U_1)=\pi_1(U_2)$, which is the fundamental group of $f_i$-fiber.

To avoid a technique complexity, we will state the most general gluing criterion in \cite{CFG}, because our gluing charts are of the simplest form, a trivial bundle with fiber a fixed infra-nilamnifold, a gluing technique presented in Appendix 2 of \cite{CFG} is
all we need in the proof of Theorem B (Method 2).
\endremark

\vskip4mm

\subhead 1.3. Cheeger-Colding-Naber theory on Ricci limit spaces
\endsubhead

\vskip4mm

A length space, $(X,p)$, is called a Ricci limit space, if there is a sequence of pointed complete $n$-manifolds, $(M_i,p_i)@>\op{GH}>>(X,p)$,
with $\op{Ric}_{M_i}\ge -(n-1)$.

Observe that the most significant difference,  from a limit space of a sequence of $n$-manifolds $M_i$
with sectional curvature, $\op{sec}_{M_i}\ge -1$, is a week link between large scale geometry and
small scale geometry; for instance, the tangent cone at $x\in X$ is not unique, and may not be a
metric one.

Two basic tools in the Cheeger-Colding-Naber theory are two quantitative rigidity results on manifolds
with Ricci curvature bounded below:  the Splitting theorem of Cheeger-Gromoll, and the rigidity maximal volume rigidity.

\proclaim{Theorem 1.3.1} {\rm (Quantitative splitting, \cite{CC1}, \cite {CoN})} Given $n, \delta, \epsilon, L>0$, there exists a constant
such that if $M$ is a complete $n$-manifold, $p, q, r\in M$ such that $E(p)=d(p,q_-)+d(p,q_+)-d(q_-,q_+)<\epsilon$,
$d(p,q_\pm)\ge L$, then
$$d_{\op{GH}}(B_R(p),B_R(X\times \Bbb R))<\Psi(\delta,\epsilon,L^{-1}|n,R),$$
where $0<R<<L$ and $X$ is a compact length space.
\endproclaim

\proclaim{Corollary 1.3.2} {\rm (Splitting of Ricci limit space)} Let $(M_i,p_i)@>\op{GH}>>(X,p)$ such that $\op{Ric}_{M_i}\ge -(n-1)\delta_i\to 0$.
If $X$ contains a line, then $X$ splits with a $\Bbb R^1$-factor, $X=Y\times \Bbb R^1$.
\endproclaim

\proclaim{Theorem 1.3.3} {\rm(Volume convergence, \cite{Co1}, \cite{CC1})} Let $(M_i,p_i)@>\op{GH}>>(X,p)$
be sequence of complete $n$-manifolds such that
$$\op{Ric}_{M_i}\ge -(n-1),\qquad \op{vol}(B_1(p_i))\ge v>0.$$
For any $M_i\ni x_i\to x\in X$ and any $R>0$,
$\op{vol}(B_R(x_i))\to \op{Haus}^n(B_R(x))$ as $i\to \infty$, where
$\op{Haus}^n$ denotes the $n$-dimensional Hausdorff measure.
\endproclaim

\proclaim{Theorem 1.3.4} {\rm(Almost volume balls are almost metric balls, \cite{CC1})} Given $n, \epsilon>0$, there exists $\delta=\delta(n,\epsilon)>0$ such that
if a complete $n$-manifold $M$ satisfies
$$\op{Ric}_M\ge -(n-1)\delta,\quad \op{vol}(B_1(x))\ge (1-\delta)\op{vol}(\b B_1^n(0)),$$
then $d_{\op{GH}}(B_r(x),\b B_r^n(0))<\epsilon\cdot r$ for all $0<r\le 1$ (we refer $x$ as
and a $(n,\epsilon,1)$-Reifenberg point).
\endproclaim

For any $q\in X$ and any $\epsilon_i\to 0$, passing to a subsequence we call the limit space,
$(\epsilon_i^{-1}X,q)@>\op{GH}>>(T_qX,q)$ a tangent cone of $X$ at $q$. A tangent cone may depend on
the subsequence. A point $q$ is called regular, if $T_qX$ is unique and is isometric to $\Bbb R^k$, $k\le n$. We call $\Cal S=X-\Cal R$ the singular set.

\proclaim{Theorem 1.3.5} {\rm(Regularity/singularity)} Let $(M_i,p_i)@>\op{GH}>>(X,p)$ with $\op{Ric}_{M_i}\ge -(n-1)$. Then

\noindent {\rm (1.3.5.1) (\cite{CC1})} The set of regular points $\Cal R$ is dense in $X$.

\noindent {\rm (1.3.5.2) (\cite{CC1})} If $\op{vol}(B_1(p_i))\ge v>0$, then $\dim_H(\Cal S)\le n-2$.

\noindent {\rm (1.3.5.3) (\cite{ChN})} If $\op{vol}(B_1(p_i))\ge v>0$ and $\op{Ric}_{M_i}\le n-1$,
then $\Cal R$ is an open $C^{1,\alpha}$-manifold, and $\dim_H(\Cal S)\le n-4$.
\endproclaim

\proclaim{Theorem 1.3.6} {\rm(Diffeomorphic stability, \cite{CC2}, \cite{CoN})}  {\rm (1.3.6.1)} Let $M_i@>\op{GH}>>M_0$ such that
$M_i$ is a compact $n$-manifold with $\op{Ric}_{M_i}\ge -(n-1)$ and $M_0$ is a compact $n$-manifold. Then for $i$ large,
$M_i$ is diffeomorphic to $M_0$ by a $\epsilon_i$-GHA, $\epsilon_i\to 0$.

\noindent {\rm (1.3.6.2)} Let $(M_i,p_i)@>\op{GH}>>(M_0,p)$ such that $M_i$ is a complete $n$-manifold with $\op{Ric}_{M_i}\ge -(n-1)$ and $M_0$ is a complete $n$-manifold.  Given any $R>0$, for  $i$ large $B_R(p_i)$ is diffeomorphic to $B_R(p)$
by a pointed $\epsilon_i$-GHA.
\endproclaim

\proclaim{Theorem 1.3.7} {\rm(Limit groups are Lie groups, \cite{CC2}, \cite{CoN})} Let $X$ be a Ricci limit space. Then the isometry group of $X$ is a Lie group.
\endproclaim

We mention that in this paper,  (1.3.5.2) and (1.3.5.3) are not used.

\vskip4mm

\subhead 1.4. The Margulis Lemma
\endsubhead

\vskip4mm

In our proof of Theorem 0.3,  we will perform a successive blow-ups on a suitable bounded
covering space of $M_i$ in Theorem 0.3' so that limit spaces are flat tori; which is based
the generalized Margulis lemma in \cite{KW}, see Theorem 1.4.2.

\proclaim{Theorem 1.4.1} {\rm (Short generators estimate, \cite{KW})} Let $M$ be a complete $n$-manifold with $\op{Ric}_M\ge -(n-1)$.
For any $p\in M$, there is a subset of almost full measure, $A\subseteq B_1(p)$, such that for $q\in A$, the short generators for
$\pi_1(M,q)$ with length $\le R$ is bounded above by a constant $C(n,R)$.
\endproclaim

Theorem 1.4.1 is required in proving the following:

\proclaim{Theorem 1.4.2}  {\rm (Generalized Margulis lemma, \cite{KW})}  Given $n$, there exist constants, $\epsilon(n), c(n) >0$,
such that if $M$ is a complete $n$-manifold with $\op{Ric}_M\ge -(n-1)$, then for any $x\in M$, the subgroup of $\pi_1(B_1(x),x)$
generated by loops of length $\le \epsilon(n)$ contains a normal nilpotent subgroup of index $\le c(n)$.
\endproclaim

We now give a brief description on a general approach to a Margulis lemma (\cite{KW}, \cite{KPT}, \cite{FY}, \cite{XY}),
assuming Theorem 1.4.1.

Let's start with an algebraic criterion for a finitely generated discrete group $\Gamma$ to contain a nilpotent subgroup of a finite index:

\vskip1mm

\noindent (1.4.3.1) There is a descending sequence of normal subgroups,
$$\Gamma\triangleright\Gamma_1\triangleright \cdots \triangleright \Gamma_s\triangleright \Gamma_{s+1}=e,\quad \Gamma_j/\Gamma_{j+1}\cong \Bbb Z^{k_j},\quad 1\le j\le s.$$

\vskip1mm

\noindent  (1.4.3.2) The induced homomorphism via conjugation, $\rho_j: \Gamma\to \op{Aut}(\Gamma_j/\Gamma_{j+1})$,
has a finite image.

\vskip2mm

Let $\Gamma_p(\epsilon)$ denote the subgroup generated by loops at a point $p$ of length $<\epsilon$. By Theorem  1.4.1,
one may assume that $p$ is chosen such that the number of the Gromov short generators for $\Gamma_p(\epsilon)$
is bounded by $c(n)$.  A basic property is that $\Gamma_p(\epsilon)$ contains a normal subgroup of bounded
index such that the quotient group contains no element of uniformly small displacement on any open subset
of the quotient space, which is based on that the isometry group of a Ricci limit space is a Lie group (\cite{CC1},
\cite{CoN}).  We now use $\Gamma(\epsilon)$ to denote its normal subgroup of bounded index.

With $0<\epsilon\le \epsilon(n)$, via a successively blowing ups argument
one sees that the Gromov's short generators for $\Gamma(\epsilon)$ with norm in a decreasing order and
$|\gamma_{1,1}|\le \epsilon(n)$,  are naturally divided in terms of gaps in their length (norm),
$$\gamma_{1,1},...,\gamma_{1,k_1},\quad \gamma_{2,k_1+1},..., \gamma_{2,k_1+k_2},\quad \cdots ,\quad \gamma_{s,k_{s-1}+1},...,\gamma_{s,k_1+\cdots +k_s},$$
such that
$$1\le \frac{|\gamma_{i,\alpha+k}|}{|\gamma_{i,\alpha}|}\le c(n),\quad \frac{|\gamma_{i+k,\alpha}|}{|\gamma_{i,\beta}|}\to 0,$$
Let $\Gamma_j$ denote the subgroup generated by short generators in blocks $s, s-1,..., j$.

The main technical work in a proof of Theorem 1.4.2 is to deform a short path from a point to another point 
so that the length distortion is proportional to the distance between the two points; once this is
achieved, one is able to verify (1.4.3.1) and (1.4.3.2). Roughly, the deformation is by the gradient flows of
certain functions constructed based on the geometry of $r^{-1}M$, where $r$ is proportional to the norm of short generators in some block.
Precisely, combining the holonomy distortion estimate and the norm gaps in blocks one is able to conclude that each $\Gamma_j$ is normal in $\Gamma$, and the holonomy representation, $\phi_j: \Gamma\to
\op{Aut}(\Gamma_j/\Gamma_{j+1})$, has a (bounded) finite image. Moreover, one can show
that $\Gamma_j/\Gamma_{j+1}$ contains a free abelian group of (bounded) finite index.

Let's analyze the gradient flows on manifolds with sectional (resp. Ricci) curvature bounded below.

\vskip1mm

\noindent (1.4.4.1) Assume that $\op{sec}_M\ge -1$. The main technique, a gradient push introduced
in  \cite{KPT} (cf. \cite{XY}), allows one deform a path between any two $(m,\delta)$-strained points
in the scale of $r^{-1}M$. A point $p$ is called an $(m,\delta)$-strained on $r^{-1}M$ if there are $2m$-points, $(a_i^\pm, b_i^\pm)$,
with $d(p,\pm a_i), d(p,\pm b_i)\ge \rho>0$, $\measuredangle (a_i^+pa_i^-)>\pi-\delta$ and
$\measuredangle (a_i^\pm pb_j^\pm )>\frac \pi2-\delta$.  The gradient push is the gradient flows of the
distance function at some point from $(a_i^\pm, b_i^\pm)$, and the distance
distortion is exponentially proportional to the time of gradient flows. Note that along the gradient push one has
to change distance functions from time to time, determined by each blow-up rate, the set of $(m,\delta)$-strained point on $r^{-1}M$
is path connected, open and dense, hence one can further gradient push a path to a point which is not $(m,\delta)$-strained
(i.e., relative singular).

\vskip1mm

\noindent (1.4.4.2) Assume that $\op{Ric}_M\ge -(n-1)$. In general, there is no control from a large scale geometry on a small
scale geometry (e.g., contrasts to $\op{sec}_M\ge -1$, the sizes of a geodesic triangle do not control angles), so
the gradient-push technique is replaced by the so-called a zooming-in map in \cite{KW}, based on $(k,\delta)$-splitting
maps (\cite{CC}). With a similar blow-up rate as in (1.4.1.1), one obtains a $(k,\delta)$-splitting map, whose gradients are almost orthonormal
in the $L^2$-sense.  The weak control on gradient leads to a possible vanishing of a gradient at some point, or a $(k,\delta)$-splitting map on $B_1(p)$ may not be a $(k,\delta)$-splitting map at $q\in B_1(p)$ in a smaller scale at $q$. As a result, for the holonomy
estimate to hold, one has to exclude certain points.

\remark{Remark \rm 1.4.5} In \cite{CJN}, it was proved that an $(n,\delta)$-splitting map is
everywhere non-degenerate. In our proof of Theorem A, we show that the points on the local Riemannian universal
cover of any $\rho$-ball in $M_i$ are uniform $(n,\delta)$-Reifenberg (Lemma 2.3), which is used to
check that an $(k,\delta)$-splitting map is non-generate everywhere (\cite{Hu}). Combining the transformation
technique in \cite{CJN}, it is expected that the zooming in maps works at every point, hence
a holonomy estimate should hold at every point, similar to (1.4.4.1). Because our proof of Theorem 0.3' does not require this property,  we will not pursuit a proof in this paper.
\endremark

\vskip4mm

\subhead 1.5. Collapsing with/without local bounded Ricci covering geometry
\endsubhead

\vskip4mm

We will supply examples mentioned in the previous discussions.

\example{Example 1.5.1} {\rm (No nearby metrics with bounded sectional curvature)}
Given any sequence, $M_i@>\op{GH}>>X$,
of compact $n$-manifolds such that
$$\op{Ric}_{M_i}\ge -(n-1),\quad \op{diam}(M_i)\le d, \quad \op{vol}(M_i)\ge v>0.$$
Then the metric product, $M_i\times \epsilon_i S^1@>\op{GH}>>X$ ($\epsilon_i\to 0$) satisfies the local
$(1,v)$-bounded Ricci covering geometry.

If $X$ is not a topological manifold (\cite{An}, \cite{Ot}, \cite{Per}), then
$M_i\times \epsilon_i S^1$ admits no nearby metric with bounded sectional curvature. Otherwise,
the collapsing structure on
$M_i\times \epsilon_i S^1$ with respect to the nearby metric with bounded sectional curvature has
to be the free $S^1$-action (Theorem 0.2), thus the orbit space, $(M_i\times \epsilon_i S^1)/S^1=M_i$, is homeomorphic to $X$;
a contradiction.
\endexample

\example{Example 1.5.2} ({Collapsing with smooth limits, no fibration, \cite{An})}  If one removes the volume condition in
Theorem 0.4, then there is no fiber bundle map.

Anderson (\cite{An}) constructed one family of metrics $g_\delta$ on
a compact $n$-manifold $M$ with a nontrivial second homotopy group such that
$$(M,g_\delta)@>\op{GH}>>T^{n-1} \, (\text{flat}),\quad 0<\delta\le \op{Ric}_{g_\delta}\le 1.$$
There is no $S^1$-fiber from $(M,g_\delta)$ to $T^{n-1}$; otherwise all $\pi_q(M)=0$ for $q\ge 2$.
\endexample

\example{Example 1.5.3} (Non-smooth/smooth limit vs. infra-nil/nilpotent fundamental groups)

Let $(K^2,g)$ denote a Klein bottle, and let $S^1$ acts on $K^2$ by isometries with two exceptional $S^1$-orbits with
isotropy $\Bbb Z_2$. Hence, $K^2/S^1=[0,1]$. Decomposes $g_0$ at each point as $g=ds^2+g^\perp$, where $ds^2$ denote the restriction of $g$ to the tangent space of $S^1$-orbit, and $g^\perp$ the orthogonal complement. For $\epsilon\in (0,1]$, let $g_\epsilon=\epsilon^2ds^2+g^\perp$.
Then $(K^2,g_\epsilon)@>\op{GH}>>[0,1]$ with $\op{sec}_{g_\epsilon}\equiv 0$, and $(K^2,\sqrt{\epsilon}g_\epsilon)@>\op{GH}>>\op{pt}$
with $\op{sec}_{\sqrt{\epsilon}g_\epsilon}\equiv 0$.  When blowing up by $\op{diam}(K^2,\sqrt\epsilon g_\epsilon)^{-1}\approx \epsilon^{-\frac 12}$, $(K^2,g_\epsilon)@>\op{GH}>>[0,1]$, while the double cover, $(T^2,g_\epsilon)@>\op{GH}>>S^1$.
\endexample

\example{Example 1.5.4} {\rm (Collapsing Calabi-Yau metrics along a nilpotent fibration
with singular fibers not topological manifolds)}. Let $M$ be a compact K\"ahler manifold. Assume there is
a map, $f: M\to X$, where $X$ is a smooth manifold except a finite isolated singular points,  $\{x_i\}$,
and $f: M-\bigcup_if^{-1}(x_i)\to X-\{x_i\}$ is smooth nilpotent fibration, $f^{-1}(x_i)$ is a not submanifold
of dimension $\dim(M)-\dim(X)$.
In \cite{GW}, \cite{GTZ1,2},  \cite{HSVZ}, and \cite{Zh}, a one-parameter family of Calabi-Yau metrics
$g_\epsilon$ are constructed so that $(M,g_\epsilon)@>\op{GH}>>X$. It can be checked that away from
the singular fibers, sectional curvature is bounded,  thus a local Ricci bounded covering geometry
is satisfied.

Indeed, with some minor work, a local version of Theorem B apply to the above situation,
$f: M-\bigcup_if^{-1}(x_i)\to X-\{x_i\}$ is a nilpotent fiber bundle.
\endexample

\vskip4mm

\head 2. Outline of Proof of Theorem A
\endhead

\vskip4mm

In this section, we will outline our proof of Theorem A, with details supplied in Section 3.

By the Gromov's pre-compactness, Theorem 0.3 has  the following equivalent form (which we will
present a proof independent of Theorem 0.1):

\proclaim{Theorem 0.3'} Let a sequence of compact $n$-manifolds, $M_i@>\op{GH}>>\op{pt}$, such that
$$\op{Ric}_{M_i}\ge -(n-1), \quad \op{vol}(B_1(\tilde p_i))\ge v>0, \quad \tilde p_i\in \tilde M_i.$$
Then for $i$ large, $M_i$ is diffeomorphic to an infra-nilmanifold.
\endproclaim

A starting point is Theorem 1.4.2, by which we may assume that a bounded normal covering of $M_i$, $M_{i,0}$, with $\pi_1(M_{i,0})$
a nilpotent group.  The main effort in our proof of Theorem 0.3' is to show that $M_{i,0}$ is
diffeomorphic to a nilmanifold,  which replies on the following fiber bundle criterion for a nilmanifold.

\proclaim{Theorem 2.1} {\rm (Nilmanifolds, a sequence of bundles)}  A compact manifold $M_0$  is
diffeomorphic to a nilmanifold  if and only if $M_0$ carries a sequence of bundles:
$$M_1\to M_0\to B_1,\quad M_2\to M_0\to B_2,\quad \cdots, \quad M_{s+1}\to M_0
\to B_{s+1},$$
such that

\noindent {\rm (2.1.1)} Each $M_i$-fiber is a union of $M_j$-fibers for all $j\ge i$, $M_{s+1}=\{\op{pt}\}$.

\noindent {\rm (2.1.2)} For each $1\le i\le s$, the natural projection, $B_{i+1}\to B_i$, is a principal $T^{k_i}$-bundle.
\endproclaim

We will show that $M_0$ in Theorem 2.1 admits an iterated principal circle bundles:
$$S^1\to M_0\to M_1,\quad S^1\to M_1\to M_2,\quad \cdots, \quad  S^1\to M_n\to \op{pt},$$
which, according to \cite{Na} and \cite{Be} (Theorem 3.1), implies that $M_0$ is diffeomorphic to
a compact nilmanifold.

Our main work is the following result.

\proclaim{Theorem 2.2}  Let the assumptions be as in Theorem 0.3', and let $M_{i,0}$ be a (bounded) finite normal covering space of $M_i$ such that $\pi_1(M_{i,0})$ is nilpotent.
Then for $i$ large $M_{i,0}$ admits a sequence of bundles satisfying  Theorem 2.1 i.e. $M_{i,0}$ is diffeomorphic to a nilmanifold.
\endproclaim

A by-product of the proof of Theorem 2.2 is that we are able define, on $M_{i,0}$, a left-invariant metric, such that the normal covering transformations on $M_{i,0}$ by
$\Gamma_i/\Gamma_{i,0}$ are almost isometries. Using the left invariant metric, we apply the center of mass technique (Lemma 3.5.1) to modify the normal covering transformations so that the $\Gamma_i/\Gamma_{i,0}$ action preserves the nilpotent structure on $M_{i,0}$ defined by the sub-bundles in Theorem 2.2.

As a preparation, we need some lemmas.

\proclaim{Lemma 2.3} {\rm (Euclidean limits and uniform Reifenberg points)}
Let $M_i@>\op{GH}>>\op{pt}$ with local $(\rho,v)$-bound  Ricci covering geometry. Then up to a scaling of $M_i$, passing to
a subsequence and for any $R>1$ the following hold:
$$\CD (\tilde M_i,\tilde p_i,\Gamma_i)@>\op{eqGH}>>(\Bbb R^n,0,\Bbb R^n)\\
@V \pi_iVV @V\op{proj}VV\\
M_i@>\op{GH}>>\op{pt},\endCD,\qquad \cases \op{Ric}_{M_i}\ge -(n-1)\delta_i & \delta_i\to 0\\ d_{\op{GH}}(B_{i,s}(\tilde x_i),\b B^n_s(0))<s\epsilon_i& \epsilon_i\to 0\\
\tilde x_i\in B_R(\tilde p_i) &0<s\le R\endcases.$$
\endproclaim

We will refer a maximal collapsing sequence in Lemma 2.3 as satisfying a uniform Reifenberg condition. From now on, we will always assume a sequence in Theorem 0.3' satisfying Lemma 2.3. Consequently, given any sequence, $r_i\to \infty$, and at any $\tilde x_i\in \tilde M_i$, the $r_i$-blow up limit is always $\Bbb R^n$ i.e., $(r_i\tilde M_i,\tilde x_i)@>\op{GH}>>(\Bbb R^n,0)$.

We find that if $\pi_1(M_i)$ is nilpotent, then any blow-up limit of $M_i$ must also be a
product of a flat torus and an Euclidean space. This discovery leads to the following key geometric ingredients of Theorem 2.2.

\proclaim{Theorem 2.4} {\rm (Uniform local elementary T-structures)}  Let a sequence,
$M_i@>\op{GH}>>\op{pt}$, be as in Lemma 2.3. Assume that $\Gamma_i=\pi_1(M_i)$ is nilpotent. Then passing to a subsequence, the following properties hold:

\noindent {\rm (2.4.1)} There are integers,  $(k_1,...,k_s)$, $m_j=\sum_{l=1}^jk_l$, $m_s=n$, such that for any $x_i\in M_i$, a short generators (ordered in lengths decreasing)
for $\Gamma_i$ at $x_i$, divides into groups,
$$\gamma_{i,1}(x_i),...,\gamma_{i,m_1}(x_i), \gamma_{i,m_1+1}(x_i),...,
\gamma_{i,m_2}(x_i),...,\gamma_{i,m_{s-1}+1}(x_i),...,\gamma_{i,m_s}(x_i),$$
$$\frac{l_j(x_i)}{\ell_j(x_i)}\le 1,\quad \frac {\ell_j(x_i)}{\ell_{j'}(x_i)}\to \infty, \quad j'>j,\quad 1\le j\le s,\quad (i\to \infty),$$
where
$$\cases \ell_j(x_i)=\max\{|\gamma_{i,k}(x_i)|\}_{k=m_j+1}^{m_{j+1}}\\
l_j(x_i)=\min\{|\gamma_{i,k}(x_i)|\}_{k=m_j+1}^{m_{j+1}}. \endcases$$

\noindent {\rm (2.4.2)} There exists a constant $c=c(\{M_i\})>0$ such that the following properties hold:
for any $\epsilon>0$, there is $N(\epsilon)>0$, such that for $i\ge N(\epsilon)$ and any $x_i\in M_i$,
$$d_{\op{GH}}(B_{30\sqrt n}(x_i,\ell_j(x_i)^{-1}M_{i,0}), B_{30\sqrt n}(T^{k_j}\times \Bbb R^{m_{j-1}}))<\epsilon,$$
where $T^{k_j}\times \Bbb  R^{m_j}$ is flat with $T^{k_j}$ satisfying
$$\op{diam}(T^{k_j})\le \sqrt k_j\le \sqrt n,\quad \op{injrad}(T^{k_j})\ge c$$
\endproclaim




Based on Theorem 2.4, we able to construct an atlas of elementary $T$-structures on suitable spaces, to which we apply Theorem 1.2.3 to construct, by induction on $s$, a sequence of sub-bundles on $M_{i,0}$ satisfying Theorem 2.1. Here is an outline of the proof:

Fixing $\epsilon>0$ and $i>N(\epsilon)$ in Theorem 2.4, we construct a locally bounded
open cover, $\{U_{i,\alpha}\}$, $$B_{10\sqrt n}(x_{i,\alpha},\ell_s(x_{i,\alpha})^{-1}M_{i,0})\subset U_{i,\alpha}\subset B_{11\sqrt n}
(x_{i,\alpha},\ell_s(x_{i,\alpha})^{-1}M_{i,0}), \tag 2.5.1$$
such that for any $x_i\in M_i$, $\{x_{i,\alpha}\}\cap B_{\ell_s(x_i)30\sqrt n}(x_i)$ is almost $\ell_s(x_i)10\sqrt n$-separated and
almost $\ell_s(x_i)10\sqrt n$-dense.

\noindent (2.5.2) There is a diffeomorphism and $\epsilon$-GHA (Theorem 1.3.6),
$f_{i,\alpha}: U_{i,\alpha}\to T^{k_s}_\alpha\times \b B^{m_{s-1}}_{10\sqrt n}(0_\alpha)$, where $U_{i,\alpha}$ (resp. $T^{k_s}\times \b B^{m_{s-1}}_{10\sqrt n}(0_\alpha)$) is equipped with the induced metric on $\ell_s(x_{i,\alpha})^{-1}M_{i,0}$ (resp. the product of flat metrics).
Moreover, a choice of a base point $e\in T^{k_s}$ and a canonical basis for $\pi_1(T^{k_s},e)$ uniquely determines (via parallel translation) a Lie group structure on $T^{k_s}$,  thus defines a free $T^{k_s}$-action on $U_{i,\alpha}$ with orbits $f_{i,\alpha}^{-1}(T^{k_s}_\alpha\times v)$, for any $v\in \b B_{10\sqrt n}^{m_{s-1}}(0_\alpha)$.

\noindent (2.5.3) If $U_{i,\alpha}\cap U_{i,\beta}\ne \emptyset$, where $B_{\frac 1{10}}(z_i)\subset U_{i,\alpha}\cap U_{i,\beta}$ for some $z_i$.

\noindent (2.5.4) There is a diffeomorphic $\epsilon$-GHA, $\phi_{i,\alpha,\beta}: U_{i,\alpha,\beta}\to T^{k_s}\times \b B_{30\sqrt n}^{m_{s-1}}(0)$, such that $\psi_{i,\alpha}=\phi_{i,\alpha,\beta}\circ f_{i,\alpha}^{-1}: T^{k_s}_\alpha\times
\b B_{10\sqrt n}^{m_{s-1}}(0)\to T^{k_s}\times \b B_{30\sqrt n}^{m_{s-1}}(0)$ is an diffeomorphic embedding and $2\epsilon$-GHA.

It is obvious that the two $T^{k_s}$-actions on $U_{i,\alpha}\cap U_{i,\beta}\ne \emptyset$ are $C^0$-close (with respects to
either the original or the pull-back flat metric). Because the free $T^{k_s}$-action on $U_{i,\alpha}$ are isometric with
respect to the pull-back metric from $T^{k_s}\times \Bbb R^{m_{s-1}}$, the $C^1$-close condition in
Theorem 1.2.3 is satisfied ( Remark 1.2.4).
Hence, we are able to glue the atlas of elementary $T^{k_s}$-structures, $\{(U_{i,\alpha},T^{k_s})\}$, to form
an affine $T^{k_s}$-bundle on $M_{i,0}$, $T^{k_s}\to M_{i,0}@>h_{s,i}>>B_{i,s}$. Because $\Gamma_{i,0}$ is torsion-free nilpotent and
$\Gamma_{i,s}$ is normal abelian subgroup, $\Gamma_{i,s}$ is contained in the center of $\Gamma_{i,0}$. Hence, the holonomy group
is trivial i.e., the $T^{k_s}$-bundle is principal.

Next, we shall construct a principal $T^{k_{s-1}}$-bundle on $B_{i,s}=M_{i,0}/T^{k_s}$, where the free $T^{k_s}$-action on $M_{i,0}$
are almost isometries with respect to $\ell_s(x_{i,\alpha})^{-1}g_i$ around $x_{i,\alpha}$ (hence with respect to any $\ell_j(x_{i,\alpha})^{-1}g_i$, $j=1,...,s-1$).

Fixing a probability measure on $T^{k_s}$, let $g_{i,s}$ denote the $T^{k_s}$-invariant metric obtained by averaging $g_i$ under the $T^{k_s}$-action. Then for any fixed $R>1$,
$$d_{\op{GH}}(B_R(x_i,\ell_s(x_i)^{-1}g_{i,s}), B_R(x_i,\ell_s(x_i)^{-1}g_i))\to 0,\quad i\to \infty,$$
because $\frac{\ell_s(x_i)}{\ell_{s-1}(x_i)}\to 0$, all $T^{k_s}$-orbits on the above $R$-balls (with respect to both metrics)
uniformly collapse to points.  Hence, the above implies the following,
$$d_{\op{GH}}(B_R(\bar x_i,\ell_{s-1}(x_i)^{-1}\bar g_{i,s}), B_R(T^{k_{s-1}}\times \Bbb R^{m_{s-2}}))\to 0,\quad i\to \infty.$$

Consider the blow-up at the second fastest rate
at any $x_{i,\beta}$ , $\ell_{s-1}(x_{i,\beta})^{-1}M_{i,0}$, again by Theorem 2.4 we can construct a locally finite open cover for $M_{i,0}$,
$\{V_{i,\beta}\}$, where similar to (2.5.1) $V_{i,\beta}$ is $T^{k_s}$-invariant open subset of $M_{i,0}$,
$$B_{10\sqrt n}(x_{i,\beta},\ell_{s-1}(x_{i,\beta})^{-1}M_{i,0})\subset V_{i,\beta}\subset
B_{11\sqrt n}(x_{i,\beta},\ell_{s-1}(x_{i,\beta})^{-1}M_{i,0}).$$

\proclaim{Lemma 2.6} Let $B_{i,s}=M_{i,0}/T^{k_s}$ be equipped with $\bar g_{i,s}$. Then $B_{i,s}$ admits a locally bounded open cover, $\{\bar V_{i,\beta}\}$, and a diffeomorphic $\epsilon_i$-GHA ($\epsilon_i\to 0$ as $i\to \infty$),
$$f_{i,\beta}:  \ell_{s-1}(x_{i,\beta})^{-1}\bar V_{i,\beta}\to T^{k_{s-1}}\times B_{10\sqrt n}(\Bbb R^{m_{s-2}}),$$
and $T^{k_{s-1}}\times B_{10\sqrt n}(\Bbb R^{m_{s-2}})\subset T^{k_{s-1}}\times \Bbb R^{m_{s-2}}$ is a metric product of flat manifolds.
\endproclaim

Consequently, we obtain
a locally bounded open cover for $B_{i,s}$, $\{\bar V_{i,\beta}\}$, satisfying properties as in (2.5.2)-(2.5.4).
By now we are able to apply Theorem 1.2.3, gluing the local elementary $T$-structures, $\{(\bar V_{i,\beta}, T^{k_{s-1}})\}$,
on $B_{i,s}$, and obtain a principal $T^{k_s}$-bundle, $T^{k_s}\to B_{i,s}@> h_{i,s-1}>>B_{i,s-1}$, and a fiber bundle,  $M_{i,s}\to M_{i,0}@>h_{i,s-1}\circ h_{s,i}>>B_{i,s-1}$, and each $M_{i,s}$ is a union of $T^{k_s}$-orbits.

Observe that if $s=2$, then the above prove Theorem 2.2. For $s>2$, we will proceed with induction on $s$ to finish
the proof of Theorem 2.2.

\remark{Remark \rm 2.7} (2.7.1) Theorem 2.4 (both (2.4.1) and (2.4.2)) will be false, if one removes the assumption that
$\pi_1(M_i)$ is nilpotent (Example 1.5.3).

\noindent (2.7.2) When restricting to the case with $|\op{sec}_{M_i}|\le 1$,  by Theorem 0.2 one sees a uniform bound on
$\frac{\ell_j(x_i)}{\ell_j(y_i)}$ for any $x_i, y_i\in M_i$. Note that our proof of Theorem 2.2 does not required such a (uniform) bound, because our proof of Theorem 2.2 relies on local
geometry (comparing the discussion at the end of Section 1.4).

\noindent (2.7.3) Let $M_i$ be as in Theorem 0.3'. Combining Theorem 1.4.2 and (2.4.1), one sees that for any $x_i\in M_i$,
in Theorem 0.3', the number of short generators at $x_i$ is uniformly bounded above by a constant depending on $n$
(comparing to Theorem 1.4.1 where the bound only holds at points which may not be specified).
\endremark

\vskip4mm

\head 3. Proof of Theorem A
\endhead

\vskip4mm

In this section, we will prove theorems/lemmas stated in Section 2.

\vskip4mm

\subhead 3.1. Proof of Theorem 2.1
\endsubhead

\vskip4mm

Our proof of Theorem 2.1 uses the following criterion of a nilmanifold.

\proclaim{Theorem 3.1}{\rm (Nilmanifolds: iterated principal circle bundles, \cite{Na}, \cite{Be})} A compact manifold $M$ is
diffeomorphic to nilmanifold $N/\Gamma$ if and only if $M$ admits an iterated principal circle bundles:
$$S^1\to M\to M_1,\quad S^1\to M_1\to M_2,\quad \cdots \quad S^1\to M_n\to \op{pt}.$$
\endproclaim

\demo{Proof of Theorem 2.1}

Given a nilmanifold, $N/\Gamma$, there is the descending sequence,
$$N\triangleright N_1\triangleright \cdots \triangleright N_s\triangleright N_{s+1}=e,\quad \Gamma>\Gamma_1>\cdots >\Gamma_s>\Gamma_{s+1}=e,$$
where $N_1=[N,N]$, $N_i=[N,N_{i-1}]$, $i\ge 2$, and $\Gamma_i=N_i\cap \Gamma$, and $N_s\subseteq Z(N)$, the center of $N$.
Then one obtains the desired iterated normal bundle over tori with trivial holonomy representation,
$$N_1/\Gamma_1\to N/\Gamma\to T^{k_1}=(N/N_1)/(\Gamma/\Gamma_1),\quad \cdots, \quad \op{e}\to N_s/\Gamma_s\to T^{k_s}.$$
Let $M=N/\Gamma$, $M_i=N_i/\Gamma_i$, $i\ge 1$.
Then the correspond sequence of sub-bundles are:
$$M_1\to M\to B_1=T^{k_1},\quad M_2\to M\to B_2, \quad \cdots,
\quad \op{pt}\to M\to B_{s+1}=M,$$
where $B_j=(N/N_j)/(\Gamma/\Gamma_j)$, where $T^{k_j}\to B_j\to B_{j-1}$ is a principal bundle.

We now assume that $M$ satisfies the conditions of Theorem 2.1. By Theorem 3.1, it suffices
to show that $M$ admits an iterated principal circle bundles. We proceed by
induction on $s$, starting with $s=1$, $T^{k_2}\to M\to T^{k_1}$. Because $T^{k_2}\to M\to T^{k_1}$ is
a principal $T^{k_2}$-bundle, $M$ can be expressed as iterated principal circle bundles.

Consider an iterated bundles over tori of length $s$ on $M$ in Theorem 3.1. Note that
$M_s\to M\to B_s$ is a principal $T^{k_s}$-bundle. Consider
$$T^{k_s}\to M\to M/T^{k_s}.$$
It is easy to see that $M/T^{k_s}$ again satisfies the conditions of Theorem 2.1 of length $s-1$,
and by induction $M/T^{k_s}$ admits iterated circle bundles (Theorem 3.1). By now we
can conclude that $M$ admits iterated principal circle bundles.\qed\enddemo

\vskip4mm

\subhead 3.2. Proof of Lemma 2.3 and Theorem 2.4
\endsubhead

\vskip4mm

\demo{Proof of Lemma 2.3}

Passing to a subsequence, we may assume the following commutative diagram (Lemma 1.1.3):
$$\CD (\tilde M_i,\tilde p_i,\Gamma_i)@>\op{eqGH}>>(\tilde X,\tilde p,G)\\
@V\pi_iVV  @V\op{proj}VV\\
(M_i,p_i)@>\op{GH}>>\op{pt}=\tilde X/G.
\endCD$$
By (1.3.5.1), $\tilde X$ contains a point where the tangent cone is unique and is isometric to $\Bbb R^n$. Because $G$ acts
transitively on $\tilde X$ by isometries, every point has a unique tangent cone isometric to $\Bbb R^n$.
Assume that $r_j\to \infty$ such that $(r_j\tilde X,\tilde p)@>\op{GH}>> C_{\tilde p}\tilde X=\Bbb R^n$. Then
for any $\tilde x\in \tilde  X$, $(r_j\tilde X,\tilde x)@>\op{GH}>>\Bbb R^n$.
For  each $r_j$, we may assume $i_j$ large such that for all $i\ge i_j$, $d_{\op{GH}}((r_j\tilde M_i,\tilde p_i,\Gamma_i),(r_j\tilde X_i,\tilde p))<j^{-1}$,
and $d_{\op{GH}}(r_jM_i,\op{pt})<j^{-1}$. By Lemmas 1.1.2 and 1.1.3, a standard diagonal argument one
gets the following commutative diagram,
$$\CD (r_j\tilde M_{i_j},\tilde p_{i_j},\Gamma_{i_j})@>\op{eqGH}>>(\Bbb R^n,\tilde p,\Bbb R^n)\\
@V\pi_iVV  @V\op{proj}VV\\
(r_jM_{i_j},p_{i_j})@>\op{GH}>>\op{pt}=\Bbb R^n/\Bbb R^n.
\endCD.$$

Because $(r_j\tilde M_{i_j},\tilde p_{i_j})@>\op{GH}>>(\Bbb R^n,0)$ with $\op{Ric}_{r_j\tilde M_{i_j}}\ge -(n-1)r_j^{-2}\to 0$,
for any $\epsilon>0$, by Theorems 1.3.3 and 1.3.4, for $i$ large we may assume that $\tilde p_{i_j}$
is $(R,\epsilon)$-Reifenberg
i.e.,
$$d_{\op{GH}}(B_{i,s}(\tilde p_{i_j}),\b B^n_s(0))<\epsilon\cdot s,\quad 0<s\le R,$$
where $\b B^n_s(0)$ denotes a $s$-ball in $\Bbb R^n$. Because $\pi_i(B_1(\tilde p_{i_j}))=M_{i_j}@>\op{GH}>>\op{pt}$,
the above holds when replacing $\tilde p_{i_j}$ by any point $\tilde x_{i_j}\in B_1(\tilde p_{i_j})$, hence all points on $\tilde M_{i_j}$
are uniform $(R,\epsilon)$-Reifenberg.  Now let $\epsilon\to 0$, from the above one easily sees the desired result.
\qed\enddemo


Our proof of Theorem 2.4 is divided into several lemmas.

We will first show that if  $\pi_1(M_i)$ is nilpotent, then a blow up limit space is always a metric product of a flat torus and
an Euclidean space.

\proclaim{Lemma 3.2.1} Let $M_i@>\op{GH}>>\op{pt}$ be as in Lemma 2.3. Assume that $\pi_1(M_i)$ is nilpotent.
Then passing to a subsequence the following hold:
$$\CD(\ell_1(M_i)^{-1}\tilde M_i,\tilde p_i,\Gamma_i)@>\op{eqGH}>>(\Bbb R^n,0,\Bbb Z^{k_1}\times \Bbb R^{n-k_1})\\@V\pi_iVV @V\op{proj}VV\\
(\ell_1(M_i)^{-1}M_i,p_i)@>\op{GH}>>(T^{k_1},e),
\endCD \qquad \ell_1(M_i)=\op{diam}(M_i),\tag 3.2.1.1$$
where the limit group, $\Bbb Z^{k_1}\times \Bbb R^{n-k_1}$, acts on $\Bbb R^n$ by translations.

\noindent {\rm (3.2.1.2)} $\Gamma_i$ has a normal subgroup $\Gamma_{i,1}$, $\Gamma_{i,1}\to \Bbb R^{n-k_1}$,
such that $\Gamma_i/\Gamma_{i,1}\cong \Bbb Z^{k_1}$.

\noindent {\rm (3.2.1.3)} For $x_i\in M_i$, let $\{\gamma_{i,1,j}\}$ be a finite set of short generators for $\Gamma_{i,1}$ at $x_i$. Then
$$\ell_2(x_i)=\max_j\{d_{\ell_1(M_i)^{-1}\tilde M_i}(\tilde x_i,\gamma_{i,1,j}(\tilde x_i)\}\to 0, \quad i\to \infty.$$
\endproclaim

\demo{Proof} (3.2.1.1) By Lemma 1.14 and Lemma 3.2.1, passing to a subsequence we may start with the following equivariant GH-convergence,
$$\CD (\ell_1(M_i)^{-1}\tilde M_i,\tilde p_i,\Gamma_{i,0})@>\op{eqGH}>>(\Bbb R^n,0,G_1)\\@V \pi_i VV @ V\op{proj} VV\\
(\ell_1(M_i)^{-1}M_i,p_i)@>\op{GH}>>(X_1,p).\endCD$$
Because $\Gamma_{i,0}$ is nilpotent, the limit group $G_1$ is nilpotent. We claim that $G_1$ acts
(freely) on $\Bbb R^n$ by translations, thus $G_1=\Bbb Z^{k_1}\times \Bbb R^{n-k_1}$ and
$X_1=\Bbb R^n/G_1=T^{k_1}$.

We first show that $G_1$ acts freely. If the identity isotropy group at $0\in \Bbb R^n$, $H_0\ne e$. Because $H_0$ is
a compact subgroup of a nilpotent group $G_1$, $H_0\subset Z(G_1)$, the center of $G_1$,
and $G_1/H_0$ acts isometrically on $\Bbb R^n/H_0=\Bbb R^l\times C(S^{n-l-1}_1/H_0)$,
where $C(S^{n-l-1}_1/H_0)$ is an Euclidean cone over $S^{n-l-1}_1/H_0$.
Consequently, $\Bbb R^n/G_1=\Bbb R^l\times C(S^{n-l-1}_1/H_0)/(G_1/H_0)$ is not compact,
a contradiction. Similarly, if $\Bbb R^1<Z(G_1)$ acts not by translation, then
$\Bbb R^n/\Bbb R^1=\Bbb R^l\times Y$, and $\Bbb R^n/G'=\Bbb R^l\times Y/(G_1/\Bbb R^1)$
is not compact, a contradiction, where $Y$ is non-compact manifold of $\op{sec}_Y\ge 0$
whose unique soul is the pole point. We then repeating the above to
$\Bbb R^n/Z((G_1)_0)=\Bbb R^k$, on which $G_1/Z((G_1)_0)$ acts isometrically.

Because $G_1/Z((G_1)_0)$ is nilpotent, we repeat the above argument, and eventually
get the desired result.

\vskip1mm

(3.2.1.2) Follows from (3.2.1.1) and Lemma 1.14.

(3.2.1.3) Follows from (3.2.1.2) and Lemma 1.14.
\qed\enddemo

\proclaim{Lemma 3.2.2} Let the assumptions be as in Lemma 3.2.1. Passing to a subsequence, we have
$$\CD(\ell_2(x_i)^{-1}\tilde M_i,\tilde x_i,\Gamma_i)@>\op{eqGH}>>(\Bbb R^n,0,\Bbb Z^{k_2}\times \Bbb R^{n-k_1})\\ @V\pi_iVV @V\op{proj}VV\\
(\ell_2(x_i)^{-1}M_i,x_i)@>\op{GH}>>(T^{k_2}\times \Bbb R^{k_1},e\times 0),
\endCD \tag 3.2.2.1$$
where the limit group, $\Bbb Z^{k_2}\times \Bbb R^{n-k_1}$, acts on $\Bbb R^n$ by translations.

\noindent {\rm (3.2.2.2)} $\Gamma_{i,1}$ has a subgroup $\Gamma_{i,2}$ normal in $\Gamma_i$ such that
$\Gamma_{i,1}/\Gamma_{i,2}\cong \Bbb Z^{k_2}$.

\noindent {\rm (3.2.2.3)} For $x_i\in M_i$, let $\{\gamma_{i,2,j}\}$ be a finite set of short generators for $\Gamma_{i,2}$ at $x_i$,
and let $\ell_2(x_i)=\max\{d(\tilde x_i, \gamma_{i,2,j}(\tilde x_i))\}$. Then $\ell_2(x_i)\to 0$ as $i\to \infty$.
\endproclaim

\demo{Proof}  (3.2.2.1) Similarly, passing to a subsequence we may start with the following equivariant GH-convergence,
$$\CD (\ell_1(x_i)^{-1}\tilde M_i,\tilde x_i,\Gamma_{i,1})@>\op{eqGH}>>(\Bbb R^n,0,G_2)\\@V \pi_i VV @ V\op{proj} VV\\
(\ell_1(x_i)^{-1}M_i,x_i)@>\op{GH}>>(X_2=\Bbb R^n/G_2,x).\endCD$$
Because for any $\gamma_i\in \Gamma_i-\Gamma_{i,1}$, $d(\tilde x_i,\gamma_i(\tilde x_i))\ge \op{injrad}(T^{k_1})>0$,
$\Gamma_{i,1}@>\op{GH}>>G_2$, which acts on the $\Bbb R^{n-k_1}$-factor of $\Bbb R^n=\Bbb R^{n-k_1}\times \Bbb R^{k_1}$.
Because $\Gamma_{i,1}$ is nilpotent, $G_2$ is nilpotent.  We claim that $G_2$ acts
freely by translations on $\Bbb R^{n-k_1}$, thus $G_2=\Bbb Z^{k_2}\times \Bbb R^{n-k_1-k_2}$ and
$X_2=\Bbb R^n/G_2=T^{k_2}\times \Bbb R^{k_1}$.

Without loss of the generality, we may assume the isotropy group $H_0\ne e$, $0\in \Bbb R^n$.
Because $H_0$ is compact, $H_0\subset Z(G_2)$, the center of $G_1$, and $G_2/H_0$ acts isometrically
on $\Bbb R^n/H_0=\Bbb R^l\times C(S^{n-l-1}_1/H_0)$,
where $C(S^{n-l-1}_1/H_0)$ is an Euclidean cone over $S^{n-l-1}_1/H_0$.
Consequently, $\Bbb R^n/G_1=\Bbb R^l\times C(S^{n-l-1}_1/H_0)/(G_1/H_0)$ is not compact,
a contradiction. Similarly, if $\Bbb R^1<Z(G_2)$ acts not by translation, then
$\Bbb R^n/\Bbb R^1=\Bbb R^l\times Y$, and $\Bbb R^n/G'=\Bbb R^l\times Y/(G_2/\Bbb R^1)$
is not compact, a contradiction, where $Y$ is non-compact manifold of $\op{sec}_Y\ge 0$
whose unique soul is the pole point. We then repeating the above to
$\Bbb R^n/Z((G_2)_0)=\Bbb R^k$, on which $G_2/Z((G_2)_0)$ acts isometrically.

Because $G_2/Z((G_2)_0)$ is nilpotent, we repeat the above argument, and eventually
get the desired result.

By now, similar to the proof in Lemma 3.2.1, one gets (3.2.2.2) and (3.2.2.3).
\qed\enddemo

Note that in Lemma 3.2.1, the flat torus $T^{k_1}$ is independent of $p_i$. Consequently, for $\epsilon>0$ small,
there is $N(\epsilon)>0$ such that at every point in $M_i$, the group of short generators for $\Gamma_i$,
consisting of ones whose lengths are proportional to $\op{diam}(M_i)$, can be chosen having almost
identical lengths; which are close to  seen from parallel short basis for $\pi_1(T^{k_1})$ on $T^{k_1}$.

Note that in Lemma 3.2.2,  a priori $T^{k_2}$ may depend on $x_i$.  We will show that $k_2$ is independent of $x_i$
and every $T^{k_2}$ has a bounded geometry.

\proclaim{Lemma 3.2.3} Theorem 2.4 is true for $j=1$ i.e., $k_2$ is independent of $x_i$, $\op{diam}(T^{k_2})\le \sqrt n$,
and $\op{injrad}(T^{k_2})\ge c>0$, where $c=c(\{M_i\})$.
\endproclaim

\demo{Proof}  Arguing by contradiction, given $c_j=2^{-j}\to 0$, we will pick up a contradiction subsequence of $M_i$ as follows:

(i) Fixing $c_1$, assuming $\epsilon_1>0$ and a subsequence sequence $(M_{i_1},x_{i_1})$ such that $B_{30\sqrt n}(x_{i_1},\ell_1(x_{i_1})^{-1}M_{i_1})$ is at least $\epsilon_1$-away from any $B_{30\sqrt n}(e\times 0,T^k\times \Bbb R^{k_1})$,
$\op{injrad}(T^k)\ge c_1$. Because $T^k$ has a short basis of lengths $\le 1$, $\op{diam}(T^k)\le \sqrt k<\sqrt n$.

(ii) Fixing $c_2$, assuming that there is $\epsilon_2>0$ and a subsequence of $\{(M_{i_1},x_{i_1})\}$, $(M_{i_2},x_{i_2})$ such that $B_{30\sqrt n}(x_{i_2},\ell_1(x_{i_2})^{-1}M_{i_2})$ is at least $\epsilon_2$-away from any $B_{30\sqrt n}(e\times 0,T^k\times \Bbb R^{k_1})$, $\op{injrad}(T^k)\ge c_2$.

Iterating this process, fixing $c_j$, assuming that there is $\epsilon_j>0$ and a subsequence of $\{(M_{i_{j-1}},x_{i_{j-1}})\}$, $(M_{i_j},x_{i_j})$, such that $B_{30\sqrt n}(x_{i_j},\ell_1(x_{i_j})^{-1}M_{i_j})$ is at least $\epsilon_j$-away from any $B_{30\sqrt n}(e\times 0,T^k\times \Bbb R^{k_1})$, $\op{injrad}(T^k)\ge c_j$.  For the sake of simple notion, we will use $(M_{i_i},x_{i_i})$ to denote the diagonal subsequence
from the above construction, which servers as a contradiction sequence.

Passing to a subsequence, by the proof of Lemma 3.2.2 we obtain
$$(\ell_1(x_{i_i})^{-1}M_{i_i},x_{i_i})@>\op{GH}>>(T^k\times \Bbb R^{k_1},e\times 0).$$
Fixing $j_0$ large so that $\op{injrad}(T^k)\ge c_{j_0}>0$, for all $i_i$ large such that
$$d_{\op{GH}}(B_{30\sqrt n}(x_{i_i},\ell_1(x_{i_i})^{-1}M_{i_i}),B_{30\sqrt n}(T^k\times \Bbb R^{k_1},e\times 0))<\epsilon_{j_0},$$
we get a contradiction.


Finally we show that $k_2:=k$ in the above can be chosen independent of $x_i\in M_i$. Fixing $x_i\in M_i$, by the bounded
geometry on $T^k$, it is clear that $k_2$ can be chosen for all points in a neighborhood of $x_i$; hence the set of points with $k_2$
is open.  Similarly, the set of points with $k_2$ is also close, hence we have proved that $k_2$ is independent of points in $M_i$.
\qed\enddemo

\demo{Proof of Theorem 2.4}

Without loss of generality, we may assume that $M_i@>\op{GH}>>\op{pt}$ satisfies Lemma 2.3.
We will perform successive blow-ups on the sequence with a blow-up rate determined as in Lemmas 3.2.1-3.2.3,
and we will proceed the proof by induction on the total number of the blow-ups, which is at most $n$.

By Lemma 3.2.2-3.2.3, we may assume that (2.4.1) is true for $j=1,...,t-1$.  If (2.4.2) is false on the $t$-th blow-up sequence, following the proof of Lemma 3.2.3 we get a contradiction. Moreover, $k=k_t$ is independent of $x_i\in M_i$.
\qed\enddemo

\vskip4mm

\subhead 3.3. Proof of Lemma 2.6
\endsubhead

\vskip4mm

\demo{Proof of Lemma  2.6}

We claim that given $\epsilon>0$, there exists $N(\epsilon)>0$ such that for $i\ge N(\epsilon)$,
all $\bar y_i\in \bar V_{i,\alpha}(\bar x_i,\ell_{s-1}(x_i)^{-1}\bar g_{i,s})$ are uniform $(\epsilon,1)$-Reifenberg points,
i.e.,
$$d_{\op{GH}}(B_r(\bar y_i,\ell_{s-1}(x_{i,\alpha})^{-1}\bar g_{i,s}),\b B_r^{n-k_s}(0))<\epsilon r,\quad 0<r\le 1.$$
Assuming the claim, by the standard Reifenberg method in geometry measure theory (\cite{CC1}) we conclude
that $B_{10\sqrt n}(T^{k_{s-1}}\times \Bbb R^{m_{s-2}})$ is diffeomorphic to $\bar V_{i,\alpha}$, thus by a similar
construction of $\{U_{i,\alpha}\}$ for $M_{i,0}$, one gets a desired locally bounded open cover for $B_{i,s}$.

Arguing by contradiction, assuming that for some $\epsilon_0>0$, there is a contradicting sequence,
$$d_{\op{GH}}(B_{r_i}(\bar y_i,\ell_{s-1}(x_{i,\alpha})^{-1}\bar g_{i,s}), \b B_{r_i}^{n-k_s}(0))\ge \epsilon_0 r_i,\quad 0<r_i\le 1,$$
thus
$$d_{\op{GH}}(B_1(\bar y_i,r_i^{-1}\ell_{s-1}(x_{i,\alpha})^{-1}\bar g_{i,s}), \b B_1^{n-k_s}(0))\ge \epsilon_0. \tag 2.6.1$$
Let $\op{proj}_{i,\alpha}: \widetilde V_{i,\alpha}\to V_{i,\alpha}\to \bar V_{i,\alpha}$ denote the $\Bbb R^{k_s}_1$-fiber projection
map, where $\widetilde V_{i,\alpha}\subset \widetilde{B_\rho(x_{i,\alpha})}$, and
$\Bbb R^{k_s}_1$ denotes the universal covering group of $T^{k_s}_1$ which acts freely and isometrically
on $\tilde V_{i,\alpha}$.

Fixing $R>10$, consider the following diagram (passing to a subsequence if necessary):
$$\CD (\op{proj}_{i,\alpha}^{-1}(B_{i,1})\cap B_R(\tilde y_i),\tilde g_{i,s}),\Bbb R^{k_s}_1)@> \op{GH}>>(\b B_1^{n-k_s}(0)\times \b B_R^{k_s}(0), \Bbb R^{k_s}_2)\\
@V \op{proj}_{i,\alpha} VV @V \op{proj}VV\\
B_{i,1}=B_1(\bar y_i,r_i^{-1}\ell_{s-1}(x_{i,\alpha})^{-1}\bar g_{i,s})@> ? >>\b B_1^{n-k_s}(0),\endCD\tag 2.6.2$$
where the GH-convergence is originally with respect to $\tilde g_i$, which implies with respect to $\tilde g_{i,s}$,
because the $R^{k_{s+1}}_1$-action preserves the Euclidean metric which is $\epsilon_i$-GH close to $\tilde g_i$ with respect to
any scale (Lemma 2.3).

We claim that there is a canonical identification between groups, $\Bbb R^{k_{s+1}}_1$ and $\Bbb R^{k_{s+1}}_2$, by which
the $\epsilon_i$-GH-convergence becomes an $\epsilon_i$-equivariant $\Psi(\epsilon_i)$-GHA with respect to the (pseudo) $\Bbb R^{k_{s+1}}_i$-actions. Consequently, the eqGH-convergence in (2.6.2) descends to a GH-convergence
on $\Bbb R^{k_s}$-orbit spaces, a contradiction to (2.6.1).

Recall that the $T^{k_s}_1$-action on $V_{i,\alpha}$ is defined based on Theorem 2.4; from its proof one sees a canonical
identification of the group of elements in a short basis at $x_{i,\alpha}$ with shortest lengths are identified with a canonical
basis of the $T^{k_s}$-factor on the flat product, $\b B_{10\sqrt n}^{n-k_s}(0)\times T^{k_s}$. By the Malcev's rigidity,
this determines an identification between groups, $\Bbb R^{k_s}_1$ and $\Bbb R^{k_s}_2$.

Because the $\Bbb R^{k_s}_i$-actions preserve metrics, by which the $\Bbb R^{k_s}_i$-actions are completely determined. Hence a diffeomorphism and $\epsilon_i$-GHA is also $\Psi(\epsilon_i)$-equivariant, thus the claim has been verified.
\qed\enddemo



\vskip4mm

\subhead 3.4. $M_{i,0}$ is diffeomorphic to a nilmanifold
\endsubhead

\vskip4mm

\demo{Proof of Theorem 2.2}

Fixing small $\epsilon>0$, without loss of the generality we may assume that $M_i@>\op{GH}>>\op{pt}$ is chosen so that
Theorem 2.4 holds.

We begin by summarizing the discussion following Theorem 2.4:  with respect to the scaling functions,
$\ell_s(\cdot)^{-1}$,  $i\ge N(\epsilon)$,  by Theorem 2.4 we are able to construct, on $M_{i,0}$, a locally bounded open cover
with an elementary $T^{k_s}$-structure, $\{(U_{i,\alpha},T^{k_s}_\alpha)\}$, where $T^{k_s}$ acts freely on $U_{i,\alpha}$ and properties (2.5.1)--(2.5.4) hold.  Equipped $U_{i,\alpha}$ with the pullback flat metrics,  we are able to apply
Theorem 1.2.3 to glue together the elementary $T^{k_s}$-structures, $\{(U_{i,\alpha},T^{k_s})\}$,
to form a $T^{k_s}$-bundle, $T^{k_s}\to M_i\to B_{i,s}$, with affine structural group.  Because $\Gamma_{i,0}$ is torsion-free nilpotent
and the fiber fundamental group is a normal abelian subgroup, which has to be contained in the center of
$\Gamma_{i,0}$, hence the structural group is trivial i.e., this is a principal $T^{k_s}$-bundle,
$T^{k_s}\to M_i@>h_{i,s}>>B_{i,s}$.

Next,  with the scale function, $\ell_{s-1}(\cdot)^{-1}$, by (2.4.2) $B_{30\sqrt n}(x_{i,\alpha},\ell_{s-1}(x_{i,\alpha})^{-1}M_i)$ is close to
$B_{30\sqrt n}(T^{k_{s-1}}\times \Bbb R^{m_{s-2}})$ (note that a large radius, $30\sqrt n$, is for a gluing purpose later). By (2.4.1),
any $T^{k_s}$-fiber in $B_{30\sqrt n}(x_{i,\alpha},\ell_s(x_{i,\alpha})^{-1}M_i)$ is almost collapses to a point, as in (2.5.1)
we may assume a $T^{k_s}$-invariant open subset $V_{i,\alpha}$ $\epsilon_i$-close to $B_{10\sqrt n}(T^{k_{s-1}}\times \Bbb R^{n-k_s})$. By Lemma 2.6, similar to the construction of $\{(U_{i,\alpha},T^{k_s})\}$ for $M_{i,0}$, we
construct a locally bounded open cover for $B_{i,s}$ with elementary $T^{k_{s-1}}$-structures, $\{(\bar V_{i,\beta},T^{k_s})\}$, which satisfies properties (2.5.1)--(2.5.4). Again by Theorem 1.2.3, we are able to glue together the elementary $T^{k_{s-1}}$-structures
to form a $T^{k_{s-1}}$-bundle, $T^{k_{s-1}}\to B_{i,s}@>h_{i,s-1}>> B_{i,s-1}$, which is principal because $\pi_1(B_{i,s})$ is
torsion-free nilpotent.  Clearly, the fiber bundle, $M_{i,s-1}\to M_{i,0}@> h_{i,s-1}\circ h_{i,s}>>  B_{i,s-1}$, has the principal $T^{k_s}$-bundle, $T^{k_s}\to M_{i,0}@>h_{i,s}>>B_{i,s}$, as a sub-bundle, and each $M_{i,s-1}$ is a
union of $T^{k_s}$-fibers.

We now process by induction on $s$ (in the reverse order), assuming at up to local scales, $\ell_s(\cdot )^{-1},..., \ell_j(\cdot)^{-1}$,  we construct a sequence of bundles on $M_{i,0}$, and principal tori bundles on $B_{i,s}, B_{i,s-1},..., B_{i,j}$, such that
$$\cases T^{k_s}\to M_i@>h_{i,s}>>B_{i,s}, & \\ M_{i,s}\to M_i@>h_{i,s-1}\circ h_{i,s}>>B_{i,s-1}, & T^{k_{s-1}}_i\to B_{i,s}@>h_{i,s-1}>>
B_{i,s-1},\\
\cdots, & \cdots,\\
M_{i,j}\to M_i@>h_{i,j-1}\circ \cdots \circ h_{i,s-1}\circ h_{i,s}>>B_{i,j-1}, & T^{k_{j-1}}_i
\to B_{i,j} @> h_{i,j-1}>>B_{i,j-1},
\endcases$$
and each $M_{i,j-1}$ is a union of $M_{i,j}$-fibers.

By Theorem 2.4 with the local scale number, $\ell_{i,j-1}(\cdot)^{-1}$, we observe the following properties:
$$d_{\op{GH}}(B_{10\sqrt n}(x_i,(\ell_{i,j-1}(x_i))^{-1}M_{i,0}),B_{10\sqrt n}(e\times 0, T^{k_{j-1}}_{x_i}\times \Bbb R^{m_{j-2}}))< \epsilon.$$
Now we are in a situation similar to the scale $\ell_s(\cdot)^{-1}$, where the fiber of the bundle, $M_{i,j}\to M_{i,0}\to B_{i,j}$,
is collapsed to a point with respect to $\ell_{j-1}(\cdot)^{-1}$, thus (following the proof of Lemma 2.6) we may identify $B_{10\sqrt n}(T^{k_{j-1}}\times \Bbb R^{m_{j-2}})$ with an open subset of $B_{i,j}$, $h_{j,s}(W_{i,\gamma})=\bar W_{i,\gamma}$, equipped with a flat metric.
Consequently, we can construct a locally bounded open cover for $B_{i,j-1}$ with elementary
$T^{k_{j-1}}$-structures, $\{(\bar W_{i,\gamma},T^{k_{j-1}})\}$, which satisfies properties
as in (2.5.1)-(2.5.4). Then by applying Theorem 1.2.3 we construct the desired bundle map,
$T^{k_{j-1}}\to B_{i,j-1}@>h_{i,j-2}>>B_{i,j-2}$, which is
a principal $T^{k_{j-1}}$-bundle because $\pi_1(B_j)$ is torsion-free nilpotent. Moreover, we extend the sequence of bundles,
$M_{i,j-1}\to M_{i,0}@>h_{i,j-2}\circ \cdots \circ h_{i,s}>>B_{i,j-2}$.  By the induction, the proof is complete.
\qed\enddemo

\vskip4mm

\subhead 3.5. Completion of Proof of Theorem 0.3'
\endsubhead

\vskip4mm

Let $M_i@>\op{GH}>>\op{pt}$ be as in Theorem 0.3'. By Theorem 2.2, for $i$ large, a bounded
normal covering space, $M_{i,0}$, is diffeomorphic to a nilmanifold. To prove that $M_i$
is diffeomorphic to an infra-nilmanifold, we shall use the stability of isometric Lie group
actions in \cite{Hu}.

Consider a compact Lie group $G$ acting effectively on a manifold $M$, if two $G$-actions
on $M$ are $C^1$-close (\cite{Pa}), then the two $G$-actions are conjugate via a diffeomorphism that is isotopy to $\op{id}_M$. According to \cite{GK}, if two $G$-actions are isometric with respect to two metrics, $g_0, g_1$ on $M$, with $|\op{sec}_{g_i}|\le 1$ and the convexity radius of $g_0$ $\ge i_0>0$, then a $C^0$-close of the two $G$-actions implies the $C^1$-closeness. Moreover, one obtains a conjugation close to $\op{id}_M$ constructed via the center of mass technique. The above has been extended to a non-compact nilpotent Lie group $G$ in \cite{CFG}
(see discussion following Corollary 1.2.5.)

According to \cite{Hu}, the above stability result in \cite{GK} holds if the condition, `$|\op{sec}_{g_1}|\le 1$', is weakened to that $\op{Ric}_{g_1}\ge -(n-1)$ and all points
in $(M,g_1)$ are local rewinding $(\rho,\epsilon)$-Reifenberg points.

\proclaim{Lemma 3.5.1} Let $M$ be an $n$-manifold equipped with two complete metrics, $g_1$ and $g_0$, satisfying
$$\op{Ric}_{g_1}\ge -(n-1), \quad |\op{sec}_{g_0}|\le 1, \quad \op{injrad}(g_0)\ge i_0>0.$$
Then there exists a constant, $\epsilon(n,i_0)>0$, such that if two isometric Lie group $G$-actions with respect to $g$ and $g_0$ respectively, $\mu_1$ and $\mu_0$, are $\epsilon$-close i.e.,
$$\sup\{d_{g_0}(\mu_1(g^{-1},\mu_0(g,x)),x),\,\, x\in M, \, g\in G\}<\epsilon, \quad 0<\epsilon\le \epsilon(n,i_0),$$
the two $G$-actions are conjugate via a diffeomorphism
that is $\Psi(\epsilon |n,i_0)$-isotopic to $\op{id}_M$.
\endproclaim

Lemma 3.5.1 was stated in \cite{Hu} for compact $G$; it still holds without compactness of $G$ because the conjugate map is constructed based on local geometry. Note that if $G$ is not compact, then $\mu_1$ and $\mu_0$ are $\epsilon$-close implies that the subgroup of $G$ consisting of elements such that $\mu_1(g,\cdot)=\mu_0(g,\cdot)$ has a compact quotient.

\demo{Proof of Theorem 0.3'}

Consider $M_i@>\op{GH}>>\op{pt}$ as in Theorem 0.3', $\Gamma_i=\pi_1(M_i)$. By Theorem 1.4.1, $\Gamma_i$ contains a
nilpotent subgroup $\Gamma_{i,0}$ of index $\le c(n)$. Let $M_{i,0}=\tilde M_i/\Gamma_{i,0}$. By Theorem 2.2, for $i$ large
$M_{i,0}$ admits a sequence of sub-bundles satisfying Theorem 2.1, thus there is a diffeomorphism, $f_i: N/\Gamma_{i,0}\to \tilde M/\Gamma_{i,0}$, preserving the sequence of sub-bundles, and the lifting map is an isomorphism, $\tilde f_i: N\to \tilde M_i$, of nilpotent Lie groups (determined by the sequence of sub-bundles).

By identifying $\tilde M_i$ with $N$ via $\tilde f_i$ with the $\tilde f_i^*$-pullback metric, we may assume two $\Gamma_i$-actions $\mu_i$ on $N$, the other is obtained via the Malc\'ev rigidity, $\Gamma_i<N\rtimes \op{Aut}(N)$, such that $\mu_1=\mu_2$ on $\Gamma_{i,0}$, and $\mu_1$ is the deck transformation of $\Gamma_i$, and $N/\mu_2(\Gamma_i)$ is an infra-nilmanifold.
We shall apply Lemma 3.5.1 to slightly modify $\op{id}_N$ to a $\Gamma_i$-conjugation, which descends to a diffeomorphism from $M_i$ to $N/\Gamma_i$, an infra-nilmanifold.

To apply Lemma 3.5.1, a nearby `good' metric (i.e., $g_0$) is required. We now define a left-invariant metric on $N$ as follows: by Theorem 2.4 and following its proof we may assume that a short basis for $\Gamma_{i,0}$ (at $p_i$) naturally divided into groups,
$$\cases \gamma_{i,1},...,\gamma_{i,m_1},\\
\gamma_{i,m_1+1},...,\gamma_{i,m_2},\\
\cdots \\
\gamma_{i,m_{s-1}+1},...,\gamma_{i,m_s=n},\qquad \endcases \cases \ell_1 =\max\{|\gamma_{i,1}|,...,|\gamma_{i,m_1}|\}\\ \ell_2=\max\{|\gamma_{i,m_1+1}|,...,|\gamma_{i,m_2}|\}\\ \cdots \\
\ell_s=\max \{|\gamma_{i,m_{s-1}+1}|,..., |\gamma_{i,m_s}|\}\endcases.$$
such that for $1\le j\le s$, a blowup with $\ell_j^{-1}$ yields the following commutative diagram:
$$\CD (\ell_j^{-1}\tilde M_i,\tilde p_i,\Gamma_{i,0})@>\op{eqGH}>>(\Bbb R^n,0,\Bbb Z^{k_j}\times \Bbb R^{n-m_{j-1}})\\ @VV \pi_i V @VV \op{proj} V\\
(\ell_j^{-1}M_{i,0},p_{i,0})@>\op{GH}>>(T^{k_j}\times\Bbb R^{n-m_{j-1}},e\times 0).
\endCD$$

Assume that $\ell_j^{-1}\gamma_{i,m_j+1},...,\ell_j^{-1}\gamma_{i,m_{j+1}}$ respectively converge to $\gamma_{m_j+1},...,\gamma_{m_{j+1}}$, a canonical basis on flat $T^{k_j}$. Let $v_{i,k}=\gamma_{i,k}'(0)$  with $|\gamma_{i,k}'(0)|=\op{length}(\gamma_{i,k})$ (similarly, one
gets $v_{m_j+1},..., v_{m_{j+1}}$, a basis for $T_eT^{k_j}$). Then $\{v_{i,1},...,v_{i,m_s=n}\}$
form a basis for the Lie algebra $\hbar_i=T_{\tilde p_i}\tilde M_i=T_eN$. We then obtain a left invariant metric on $N$ by specifying an inner product on a basis for $T_p\tilde M_i=T_eN$:
$$\tilde g_N(v_{i,k},v_{i,l})=\cases \ell_j(p_i)^2\left<v_k,v_l\right> & m_j+1\le k,l\le m_{j+1}\\ 0 & \text{ otherwise}\endcases $$
where $\left<,\right>$ is the Euclidean inner product on flat $T^{k_j}$ (roughly, we use the canonical basis on $T^{k_j}$ to define the inner product for elements of $\hbar_i$ in the corresponding group, while elements of $\hbar_i$ in different groups are orthogonal).

From the above commutative diagram, it is clear that restricting to any $R$-ball at $e\in N$ with respect to $\tilde g_N$ or $\tilde g_i$ ($i$ large), the two metrics are GH-close (to the limit Euclidean metric on $\b B_R^n(0)$). Indeed, $\max |\op{sec}_{\tilde g_N}|\to 0$ (otherwise, normalizing $\max|\op{sec}_{\tilde g_N}|=1$, one concludes that the left-invariant metric $\tilde g_N$ is $C^\infty$-close to an Euclidean metric; a contradiction).

Let $\Gamma_i=\bigcup_{k=1}^s \alpha_k\Gamma_{i,0}$. For fixed $\alpha=\alpha_k$, $g\in N$, and any $\gamma\in \Gamma_{i,0}$,
$$\split d(\mu_2(\alpha\gamma)(\mu_1((\alpha\gamma)^{-1},g),g)&=d(\mu_2(\alpha)\mu_2(\gamma)   (\mu_1(\gamma^{-1})\mu_1(\alpha^{-1})g)),g)\\&=d(\mu_2(\alpha)\mu_1(\alpha^{-1})g,g).\endsplit$$
Then $\{\mu_2(\gamma)^{-1}\mu_1(\gamma)(g),\,\gamma\in \Gamma_{i,0}\}\subset B_\rho(g)$, $\rho=\max\{d(\mu_2(\alpha_k)\mu_1(\alpha^{-1}_k)(g,g), \, 1\le k\le s\}$.
We now specify a representative $\alpha_i$ as the projection of $e$ to the coset, $\alpha_i\Gamma_{i,0}$ with respect to the $\tilde g_N$. Because (passing to a subsequence if necessary) the two metrics are GH-close and
$\mu_2(\alpha)^{-1}\mu(\alpha) (e)$ can be very small (than the convexity radius of $\tilde g_N$), the finite set, $\{\mu_2(\alpha_k)^{-1}\mu_1(\alpha_k)(g)\}\subset N$, has a center of mass $c_g$.

By now we are able to apply Lemma 3.5.1 to $(N,\Gamma_i,\tilde g_N)$ and $(N,\Gamma_i,\tilde g_i)$, and conclude that to $M_i$ is diffeomorphic to an
infra-nilmanifold.
\qed\enddemo

\proclaim{Corollary 3.5.2}  The finite group of isometries preserving the nilpotent structure of iterated bundles,
thus preserves each sub-bundle. In particular, the structural group of each sub-bundles is affine.
\endproclaim

\remark{Remark \rm 3.5.3} Note that the left-invariant metric on $\tilde M_i$ is explicitly defined in terms of a short basis on $M_{i,0}$,
which is different from ones via smoothing method.
\endremark

\proclaim{Lemma 3.5.4}{\rm ((8.4) in \cite{Ro4}) There are constants, $\epsilon(n), \delta(n)>0$,
such that if a compact $n$-manifold $M$ satisfies
$$|\op{sec}_M|\le 1, \quad \op{diam}(M)<\epsilon<\epsilon(n),$$
then the injectivity radius of the Riemannian universal cover, $\op{injrad}(\tilde M)\ge \delta(n)>0$.
\endproclaim

\demo{Proof of Corollary 0.5}

By Lemma 3.5.4, an almost flat manifold $M$ in Theorem 0.1, when normalizing metric by $|\op{sec}_M|\le 1$,
satisfies the conditions of Theorem 0.3. Because the proofs of both Lemma 3.5.4 and Theorem 0.3 are independent of Theorem 0.1, Theorem 0.3 implies Theorem 0.1.
\qed\enddemo

Together with Theorem 0.2, Lemma 3.5.4 has the following consequence, which is also independent of Theorem 0.1.

\proclaim{Corollary 3.5.5} \rm{(\cite{CFG}) There exist constants, $\rho(n)$ and $v(n)>0$,  such that any complete $n$-manifold $M$ with
$|\op{sec}_M|\le 1$ satisfies a $(\rho(n),v(n))$-bound covering geometry.
\endproclaim

\demo{Proof}  For any $x\in M$,  without loss of generality we may assume that $\op{vol}(B_1(x))$ is small. Then by Theorem 0.2,
we may assume the following commutative diagram,
$$\CD F\to (F(B_\rho(x)),O(n))@>\tilde f>>(Y,O(n))\\
@V\op{proj}VV @V\op{proj}VV\\ B_\rho(x)@>f>>B_\rho(x),
\endCD$$
where $Y$ is a not collapsed (open) manifold. We now choose $\rho>0$ such that $\tilde f$ is a trivial fiber
bundle i.e., $F(B_\rho(x))$ is diffeomorphic to $F\times Y$.  Because $Y$ is not collapsed, by Lemma 3.5.4
we conclude that the universal cover of $F(B_\rho(x))$ is not collapsed, thus $B_\rho(\tilde x)$ is not collapsed.
\qed\enddemo

\vskip4mm

\head 4. Proof of Theorem B
\endhead

\vskip4mm

By a standard compactness, it suffices to consider a sequence, $M_i@>\op{GH}>>N$, in Theorem B, and prove that for $i$ large, there is a smooth bundle map and $\epsilon_i$-GHA ($\epsilon_i\to 0$), $f_i: M_i\to N$, with fiber infra-nilmanifold  and affine structural group.

We first observe a basic property on local fundamental groups.

\proclaim{Lemma 4.1} {\rm (The unique local fundamental group)}  Let $M_i@>\op{GH}>>N$ be as in Theorem B.
Let $\delta=\min\{\rho,2\op{injrad}(N)\}$, and $0<\epsilon\le \epsilon(n)$,
a constant in Theorem 1.4.2. For any  $x_i\in M_i$ and $i$ large, the isomorphic class,
$$\Gamma_i(\epsilon)=\op{Im}[\pi_1(B_\epsilon(x_i),x_i)\to \pi_1(B_\delta(x_i),x_i)],$$
is independent of $x_i$.
\endproclaim

\demo{Proof} Assume the following commutative diagram:
$$\CD (B_\rho(\tilde x_i),\tilde x_i, \Gamma_{x_i}(\epsilon))@>\op{eqGH}>>(B_\rho(\tilde x),\tilde x,G)\\
@V \pi_iVV @V\op{proj}VV\\(B_\rho(x_i),x_i)@>\op{GH}>>(B_\rho(x),x).\endCD$$
In general, $\Gamma_{x_i}(\epsilon)$ contains a normal subgroup $\Gamma_i(\epsilon)$,
of bounded index such that $\Gamma_i(\epsilon)\to G_0$, the identity component of $G$, where
$\Gamma_i(\epsilon)$ is generated by a subset of short generators of $\Gamma_{x_i}(\epsilon)$ that
moves every point in $B_{\frac \rho2}(\tilde x_i)$ uniformly small.

Because $B_\rho(x)$  is Riemannian and diffeomorphic to $\b B^m_\rho(0)$, $G=G_0$ which (pseudolly) acts freely, hence
$\Gamma_{x_i}(\epsilon)\equiv \Gamma_i(\epsilon)$ i.e., the isomorphic class of
$\Gamma_{x_i}(\epsilon)$ is independent of points in $B_{\frac \rho2}(x_i)$.
\qed\enddemo

Note that Lemma 4.1 is false if one replaces $N$ with a length space, or if one removes the
local $(\rho,v)$-bound Ricci covering geometry (Examples 1.5.3 and 1.5.2).

\demo{Proof of Theorem B (Method 1)}

Let $M_i@>\op{GH}>>N$ be as in Theorem B. By Theorem 0.4 for $i$ large one may assume a smooth fiber bundle map, $F_i: M_i@>f_i>> N$, which is an $\epsilon_i$-GHA.

For any $x_i\in M_i$, by performing a successively blow up at $x_i$ at rates determined by a short basis in $\Gamma_i(\epsilon)$ (Lemma 4.1), and following the proof of Theorem A (step by step), we conclude that a fiber $F_i$ is diffeomorphic to an infra-nilmanifold, and the fiber bundle has a local trivialization compatible with the nilpotent structure on $F_i$ (because the last blow-up at $x_i$ yields a metric product limit, $X_1\times \Bbb R^{n-m}$,
$m=n-\dim(F_i)$). Because the infra-nilmanifold structure on $F_i$ is determined
a short basis for $\Gamma_i(\epsilon)$ at $x_i$, the structural group can be reduced
to affine.
\qed\enddemo

Indeed, Theorem B can also be proved without Theorem 0.4: fixing large $i$ one may
construct on $M_i$ a local finite open cover by elementary N-structures such that each chart
also carries an invariant metric determined by a short basis in the fiber so that on
overlaps these invariant metrics are $C^1$-close in terms of the gluing criterion on
N-structures in \cite{CFG} (see Remark 1.2.6). Our construction of an elementary N-structure
and an invariant metrics follow the same arguments in the proof of Theorem A. Because the proof has its own interest, we will present it below.

Let's first consider a special case in Theorem B.

\proclaim{Lemma 4.2} Let $M_i@>\op{GH}>>N$ be as in Theorem B.  Assume that for $i$ large,
$\Gamma_i(\epsilon)$ is nilpotent. Then there is a smooth fiber bundle map, $f_i: M_i\to N$, with
fiber a nilmanifold and an affine structural group, and $f_i$ is $\epsilon_i$-GHA, $\epsilon_i\to 0$.
\endproclaim

\demo{Proof} We will only outline the proof because it imitates the proof of Theorem 2.2 (step-by-step) with obvious minor modifications.

Recall the proof of Theorem 2.2, given $M_i@>\op{GH}>>\op{pt}$ with $(\rho,v)$-bound Ricci covering geometry
such that $\Gamma_i=\pi_1(M_i)$ is nilpotent, we proved that, by via a successive blow-ups at rate depends
on short basis at $x_i\in M_i$, $M_i$ admits a sequence of bundle satisfying Theorem 2.1, hence $M_i$
is diffeomorphic to a nilmanifold.  The key ingredient is that each blow-ups is a metric product of flat
$T^{k_j}\times \Bbb R^{m_{j-1}}$, such that $T^{k_j}$ has a uniformly bounded geometry, which leads
to a construction of a sequence of bundles.  For instance, the fastest rate of blow
ups leads to a principal $T^{k_{s+1}}$-bundle over $M_i$, and the slowest below-up leads
to a $F_{i,1}$-bundle on $M_i$ over $T^{k_1}$.

We now consider $M_i@>\op{GH}>>N$ with $(\rho,v)$-bound Ricci covering geometry such that
$\Gamma_i(\epsilon)$ is nilpotent.  For any $x_i\in M_i$, we perform a successive blow-ups at $x_i$
with rates determined by a short basis for $\Gamma_i(\epsilon)$ at $x_i$. The fastest rate blow-up, $\ell_s(x_i)^{-1}M_i$,
$$\CD(\ell_s(x_i)^{-1}B_1(\tilde x_i),\tilde x_i, \Gamma_i(\epsilon))@>\op{eqGH}>>
(\Bbb R^n,\tilde x, \Bbb Z^{k_{s+1}})\\
@V\pi_i|_{B_1(\tilde x_i)}VV @V\op{proj}VV\\
(\ell_i(x_i)^{-1}B_1(x_i),x_i)@>\op{GH}>>(T^{k_{s+1}}\times \Bbb R^{n-k_{s+1}},e\times 0).
\endCD$$
which implies that for $i$ large, $M_i$ has locally finite an open cover by elementary $T^{k_{s+1}}$-structures, thus by Theorem 1.2.3
we construct a bundle, $T^{k_{s+1}}\to M_i\to B_{i,s}$, with an affine structural group.
Here differences from the proof of Theorem 2.2 are that this $T^{k_{s+1}}$-bundle may not be principal, and that the last blow-up yields a $T^{k_1}$-bundle, $T^{k_1}\to B_1\to N$, because the limit space of the last blow-up is $T^{k_1}\times \Bbb R^m$ (comparing to the proof of Theorem 2.2). also yields a bundle, $F_{i,1}\to M_i\to N$, with an affine structural group. Note that because $\pi_1(F_{i,1})=\Gamma_i(\epsilon)$ is nilpotent, the restriction of a sequence of fiber bundles on $M_i$ to $F_{i,1}$ satisfies Theorem 2.1 and therefore $F_{i,1}$ is diffeomorphic to an nilmanifold.
\qed\enddemo

Note that the proof of Lemma 4.2 may not carry over if $\Gamma_i(\epsilon)$ is
an infra-nil group, because it requires something like that there is a finite normal
cover $\hat M_i$ of $M_i$ such that $\hat \Gamma_i(\epsilon)$ is a normal subgroup of
$\Gamma_i(\epsilon)$. A priori, this property unlikely holds.

Observe that the above-desired property always holds over a local trivialization of the bundle with fiber infra-nilmanifold, hence we may apply the gluing construction in \cite{CFG} to obtain a global fiber bundle with an affine structural group.

By the proof of Lemma 4.2 and the last part of the proof of Theorem 2.2 (using that $M_{i,0}$ is nilmanifold to show that
$M_i$ is infra-nilmanifold), when restricting to a local trivialization
over $B_\rho(x)\subset N$,  we obtain the following elementary N-structures.

\proclaim{Lemma 4.3} {\rm (Uniform elementary N-structures)}  Let $M_i@>\op{GH}>>N$ be as in Theorem B.
For any $\delta>0$, up to a scaling we may assume that any $2$-ball on $N$ is $\delta$-close to $\b B_2^m(0)$
in $C^3$-norm.  When $\delta$ is suitably small, for  $i$ large the following hold:

\noindent {\rm (4.3.1)} For any $x\in N$, $x_i\in M_i$, such that $x_i\to x$, there is
open subset, $B_1(x_i)\subset U_i\subset B_2(x_i)$,  and a fiber bundle
map, $\phi_i: U_i\to B_1(x)$, which is an $\epsilon_i$-GHA, $\epsilon_i\to 0$.

\noindent {\rm (4.3.2)}  A $\phi_i$-fiber $F_i$ is diffeomorphic to an infra-nilmanifold,
equipped with a product metric determined by $\Gamma_i(\epsilon)=\pi_1(F_{i,1})$.

\noindent {\rm (4.3.3)} $U_i$,  equipped with the product of the invariant metric on $F_{i,1}$ in (4.3.2) and
the Euclidean metric on $B_1(x)$ pullbacked  by $\exp_x^{-1}$, satisfies the $C^1$-closeness
condition on overlaps.
\endproclaim

\demo{Proof of Theorem B (Method 2)}

Let $M_i@>\op{GH}>>N$ be as in Theorem B. By rescaling $N$ and $M_i$ with a suitable large constant,
without loss of generality we may assume that $M_i@>\op{GH}>>N$ satisfies the conditions of Lemma 4.3.

Let $\{x_\alpha\}$ be a $\frac 14$-net on $N$. Fixing $i$ large, for each $\alpha$, applying Lemma 4.3
we obtain a locally finite open cover for $M_i$ by atlas of elementary nilpotent structures,
$\{(U_{i,\alpha},\phi_{i,\alpha})\}$.  For each $\alpha$, let $V_{i,\alpha}=\phi_{i,\alpha}^{-1}(B_{\frac 14}(x_{i,\alpha}))$.
Clearly, $\{(V_{i,\alpha},\phi_{i,\alpha}|_{V_{i,\alpha}})\}$ also forms an atlas of elementary N-structures. We will
give each $V_i$ a left-invariant metric so that the $C^1$-closeness of the gluing criterion in Appendix 2 in \cite{CFG} holds (Remark 1.2.6).

In the proof of Theorem 0.3', we construct a left-invariant metric on $M_{i,0}$, completely determined by geometry from a short basis for $\Lambda_i$ at $x_i$. Using the left invariant metric, we construct a diffeomorphism, $M_i\to N_i/\Gamma$, hence a left-invariant metric on $M_i$ (the pullback metric).

Similarly we construct, on a $\phi_{i,\alpha}$-fiber $F_i$, a left-invariant metric determined by $\Gamma_i(\epsilon)$.
By the local uniform short basis property in Lemma 3.2.3 (based on the uniform Reifenberg points on Riemannian
universal cover), it is clear that fixing a left-invariant metric on $F_i$ as the above and choosing
$\delta>0$ sufficiently small, equip $V_{i,\alpha}$ a product metric with the Euclidean metric on
$\b B_{\frac 14}^m(0)$.  Then the following properties hold: if $V_{i,\alpha}\cap V_{i,\beta}\ne \emptyset$, then
$B_{\frac 1{20}}(z_i)\subset V_{i,\alpha}\cap V_{i,\beta}$ and $V_{i,\alpha}, V_{i,\beta}\subset U_{z_i}$ which
satisfies Lemma 4.3, and it is clear that $V_{i,\alpha}$ and $V_{i,\beta}$
into the elementary N-structure, $\psi: U_{z_i}$, and these product metrics are $\epsilon_i$-close,
their pullback metrics on $\widetilde{V_{i,\alpha}\cap V_{i,\beta}}$ are $C^1$-close (away from the boundary of
$\widetilde{V_{i,\alpha}\cap V_{i,\beta}}$) in the sense of Appendix 2 in \cite{CFG} (Remark 1.2.6)). Hence, we are able to glue the local trivial elementary N-structures on $M_i$ to
obtain a global pure N-structure, which is a fiber bundle with fiber infra-nilmanifold and affine structural group.
\qed\enddemo

\vskip4mm

\head 5. Proof of Theorem C
\endhead

\vskip4mm

The following result will be used in the proof of (C2); which is not required if one assumes that
a Ricci limit space is a manifold.

\proclaim{Theorem 5.1} {\rm (\cite{Wa})} Given a sequence of complete $n$-manifolds, $(M_i,p_i)@>\op{GH}>>(X,p)$ such that
$\op{Ric}_{M_i}\ge -(n-1)$, then

\noindent {\rm (5.1.1)} $X$ is semi-locally simply connected (hence the universal cover $\tilde X$ is simply connected) i.e., for any
$x\in X$, there is $r_x>0$ such that any loop at $x$ that is contained in $B_{r_x}(x)$ is homotopically trivial in $X$;

\noindent {\rm (5.1.2)} assume $R>0$ such that $\pi_1(M_i,p_i)$ is $R$-represented i.e., there is a set of generators and
relations of bounded number with length $\le R$.  Then there is an onto homomorphism, $h_i: \pi_1(M_i,p_i)\to \pi_1(X,p)$.
\endproclaim

\remark{Remark \rm 5.2} Note that (5.1.2) was originally stated for $M_i$ with $\op{diam}(M_i)\le R$. However, (5.1.2) can be seen
from the proof in \cite{Wa}.
\endremark

\demo{Proof of Theorem C}

(C1) We start with the following commutative diagram:
$$\CD (\bar B_1(\tilde x_i),\Gamma_{x_i}(\epsilon))@>\op{eqGH}>>(\bar B_1(\tilde x),G)\\@VV \pi_iV @VV \op{proj}V\\ \bar B_1(x_i)@>\op{GH}>>\bar B_1(x)=\bar B_1(\tilde x)/G.\endCD$$
By Theorem 1.4.2, $\Gamma_{x_i}(\epsilon)$ ($\epsilon<\epsilon(n)$) contains a nilpotent subgroup index $\le c(n)$, thus
$\Gamma_{x_i}(\epsilon)$ contains a torsion-free nilpotent subgroup, $\Lambda_i$, of finite index (\cite{Ra}).
We define $\op{rank}(\Gamma_{x_i}(\epsilon)):=\op{rank}(\Lambda_i)$.
Let $(\bar B_1(\tilde x_i),\Lambda_i)@>\op{eqGH}>>(\bar B_1(\tilde x), \hat G)$,
and we shall show that $\op{rank}(\Lambda_i)\le \dim(\hat G)\le \dim(G)$,  where $\hat G$ (resp. $G$)
contains a neighborhood of the identity which isomorphically embeds into a Lie group (\cite{Wa}). Without loss of generality,
we may assume that $\op{rank}(\Lambda_i)=m\ge 1$.

Let $\gamma_{i,1},...,\gamma_{i,m}$ denote a short basis for $\Lambda_i$ in an order of length increasing. We first assume that the subgroup, $\Lambda_{i,s}\cong \Bbb Z$, generated by the first $s$ elements is normal in $\Lambda_{i,s+1}$, $\left <\gamma_{i,s+1},\Gamma_{i,s}\right>\cong \Bbb Z^2$, and $\Lambda_{i,s+1}/
\Lambda_{i,s}\cong \Bbb Z$. Let $\Lambda_{i,s}\to G_s$, and it suffices to show that $\op{rank}(\Lambda_{i,s})\le \dim (G_s)$.
Because $o(\gamma_{i,1})=\infty$ and $|\gamma_{i,1}|\to 0$, $\dim(G_s(\tilde x))\ge 1$, thus $\dim(G_s)\ge 1$.
Assume that $\op{rank}(\Lambda_{i,s})\le \dim (G_s)$, and passing to a subsequence we may assume
$(\bar B_1(\tilde p_i),\Lambda_{i,s})@>\op{eqGH}>>(\bar B_1(\tilde p),G_s)$.
By Lemma 1.1.4, $(\bar B_1(\tilde p_i)/\Lambda_{i,s},\Lambda_{i,s+1}/\Lambda_{i,s})@>\op{eqGH}>>(\bar B_1(\tilde p)/G_s,H_s)$.
Because $H_s$ contains a short generator $\gamma_{i,s+1}$ such that $|\gamma_{i,s+1}|\to 0$ and
$o(\gamma_{i,s+1})=\infty$, $\dim(H_s(\bar x))\ge 1$, thus $\dim(H_s)\ge 1$. Then
$$\split \op{rank}(\Lambda_{i,s+1})&=\op{rank}(\Lambda_{i,s})+1\le \dim(G_s)+1\\&\le \dim(G_s)+\dim(H_s)=\dim(G_{s+1}).\endsplit$$

In general, we divide short generators for $\Lambda_i$ into groups satisfying the following condition:
$$\gamma_{i,1},..,\gamma_{i,k_1},\quad \gamma_{i,k_1+1},...,\gamma_{i,k_2},\quad \cdots, \quad
\gamma_{i,k_{m-1}+1},...,\gamma_{i,k_m},$$
such that $|\gamma_{i,m}|\to 0$ as $i\to \infty$, and
$$\cases \op{rank}(\Lambda_{i,1})=1 & \Lambda_{i,1}=\left <\gamma_{i,1},..,\gamma_{i,k_1}\right>\\
\op{rank}(\Lambda_{i,2})=2& \Lambda_{i,2}=\left<\gamma_{i,1},..,\gamma_{i,k_1},\gamma_{i,k_1+1},..., \gamma_{i,k_2}\right>\\
\cdots & \\
\op{rank}(\Lambda_{i,m})=m & \Lambda_{i,m}=\left<\gamma_{i,1},..,\gamma_{i,k_1},\gamma_{i,k_1+1},...,\gamma_{i,k_2},
\cdots , \gamma_{i,k_m}\right>,\endcases$$
every $\Lambda_{i,s}$ is normal in $\Lambda_i$, and $\Lambda_{i,s+1}/\Lambda_{i,s}\cong \Bbb Z$.
Passing to a subsequence we may assume that
$(\bar B_1(\tilde x_i), \Lambda_{i,j})@>i\to \infty>> (\bar B_1(\tilde x),G_j)$, $1\le j\le m$.  Then the above proof (replacing $\gamma_{s+1}$ with $\gamma_{k_s+1}$) will go through and complete the proof.

We now assume that $n-\dim(X)=\op{rank}(\Lambda_i)\le  \dim(G)=\dim(B_1(\tilde x))-\dim(B_1(x))\le n-\dim(X)$. Then $\dim(B_1(\tilde x))=n$,
thus we can apply Theorem 1.3.3, $\op{vol}(B_1(\tilde x_i))\to
\op{vol}(B_1(\tilde x))\ge v>0$.

\vskip2mm

(C2) In the commutative diagram, by Theorem 1.3.7 $\op{Isom}(\tilde X)$ is a Lie group, thus its a closed subgroup $G$
is a Lie group and $G/G_0$ is a discrete group, where $G_0$ is the identity component of $G$. By Lemma 1.1.4, for any $0<\epsilon\le \epsilon(n)$, $\Gamma_i$ contains a normal
subgroup $\Gamma_i(\epsilon)$ such that $\Gamma_i/\Gamma_i(\epsilon)\cong G/G_0$.

By Theorem 5.1, the universal cover of $X$, $\tilde X/H$, where $H$ is a normal subgroup generated by $G_0$ and isotropy subgroups (\cite{EW}), is simply connected. Again by Lemma 1.1.4, we may assume a normal subgroup, $\Gamma_i'(\epsilon)$ such that $\Gamma_i'(\epsilon)\to H$ and $\Gamma_i/\Gamma'(\epsilon)\cong G/H$ ($i$ large) satisfying the following commutative diagram:
$$\CD (\tilde M_i/\Gamma_i'(\epsilon),p_i',\Gamma_i/\Gamma_i'(\epsilon))@> \op{eqGH}>>(\tilde X/H,p',G/H)\\
@V \pi_iVV @V \op{proj}VV\\
(M_i,p_i)@>\op{GH}>>(X,p).\endCD$$
By (5.1.2), we may assume the following commutative diagram:
$$\CD \pi_1(\tilde M_i/\Gamma_i'(\epsilon),p_i')@>\psi_i>>\pi_1(\tilde X/H, p')=e\\
@V (\pi_i)_*VV @V \op{proj}_*VV\\
\pi_1(M_i,p_i)@>h_i>>\pi_1(X,p).\endCD$$
Because $h_i$ is an onto homomorphism, $h_i$ induces an onto homomorphism,
$$\frac{\pi_1(M_i)}{[\pi_1(M_i),\pi_1(M_i)]}=H_1(M_i)@>h_i>>H_1(X)=\frac{h_i(\pi_1(M_i))}{[h_i(\pi_1(M_i)),h_i(\pi_1(M_i))]}.$$
Thus we obtain the following commutative diagram,
$$\CD H_1(\tilde M_i/\Gamma_i'(\epsilon);\Bbb Z)@>\psi_i>>H_1(\tilde X/H;\Bbb Z)=e\\
@V (\pi_i)_*VV @V \op{proj}_*VV\\
H_1(M_i;\Bbb Z)@>h_i>>H_1(X;\Bbb Z),\endCD$$
Hence, $\ker (h_i)=\op{im}((\pi_i)_*(H_1(\tilde M/\Gamma_i'(\epsilon);\Bbb Z))=\Gamma_i(\epsilon)'/([\Gamma_i,\Gamma_i]\cap \Gamma_i'(\epsilon))$. Then
$$\split b_1(M_i)-b_1(X)&=\op{rank}(\Gamma_i'(\epsilon))/([\Gamma_i,\Gamma_i]\cap \Gamma_i'(\epsilon))\\ & \le \op{rank}(\Gamma_i'(\epsilon)/[\Gamma_i'(\epsilon),\Gamma_i'(\epsilon)]).\endsplit$$
Because $\Gamma_i(\epsilon)/\Gamma_i'(\epsilon)$ is bounded for all $i$, as in (5.1.1) $\op{rank}(\Gamma_i'(\epsilon)/[\Gamma_i'(\epsilon),\Gamma_i'(\epsilon)])\le \dim(G_0)$, and ``$b_1(M_i)-b_1(X)=n-m$'' implies that $\dim(G_0)\ge m-n$, thus $\dim(B_1(\tilde x))=\dim(B_1(\tilde x))+\dim(G_0)\ge n-m+m=n$. By Theorem 1.3.3, $\op{vol}
(B_1(\tilde p_i))\to \op{vol}(B_1(\tilde p))\ge v>0$.
\qed\enddemo

\demo{Proof of Theorem 0.8}

We claim that the Riemannian universal cover of $M$ is not collapsed, by
Theorem A $M$ is diffeomorphic to an infra-nilmanifold. Because $b_1(M)=n$, $M$ is diffeomorphic to a standard torus.

We verify the claim via a contradiction argument; assuming a sequence of
$M_i@>\op{GH}>>\op{pt}$ such that $\op{Ric}_{M_i}\ge -(n-1)$ and $b_1(M_i)=n$. Then local
fundamental group of a ball radius $1$ is $\pi_1(M_i)$, thus $M_i$ satisfies local
maximal rank in (C1). Therefore the Riemannian universal cover of $M_i$ is not collapsed.
\qed\enddemo

\demo{Proof of Theorem 0.10}

By (C2),  the Riemannian universal of $M$ has a Ricci bounded covering geometry and thus we apply Theorem B to conclude a fiber bundle, $F\to M\to N$, with $F$ diffeomorphic to an infra-nilmanifold of dimension $(n-m)$ and affine structure group. From the proof of (C2),
we see that $b_1(M)-b_1(N)=n-m$ implies that $\pi_1(F)\cong \Bbb Z^{n-m}$, thus $F$ is a torus.
\qed\enddemo

\vskip4mm

\noindent {\bf Acknowledgment}. The author would like to thank Shaosai Huang for some help comments on a
draft of the present paper.

\enddocument